\PassOptionsToPackage{dvipsnames}{xcolor}
\PassOptionsToPackage{hidelinks}{hyperref}

\documentclass{article}
\usepackage{orcidlink}
\let\orcidID\orcidlink
\let\runinhead\paragraph
\title{Inferring a Cell Structure on the Space of Cyclooctane Conformations}
\author{Ulrich Bauer\orcidID{0000-0002-9683-0724}, Fabian Lenzen\orcidID{0000-0001-9579-6854}}

\usepackage{graphicx,mathtools,amsthm,array,booktabs,calc,multirow,nowidow,wrapfig,microtype,caption,subcaption}
\usepackage{newtxtext,newtxmath}
\DeclareCaptionLabelSeparator{enspace}{\enspace}
\DeclareCaptionStyle{default}[justification=centering]{font=small,labelfont=bf,labelsep=enspace,justification=justified}
\DeclareCaptionStyle{subcaption}{style=default,font=small,labelfont=normalfont,labelformat=parens,justification=centering}
\let\samenumber\ContinuedFloat

\usepackage[hidelinks]{hyperref}
\usepackage[capitalize]{cleveref}
\crefformat{equation}{(#2#1#3)} 
\crefrangeformat{equation}{(#3#1#4 – #5#2#6)}
\crefmultiformat{equation}{(#2#1#3)}{ and~(#2#1#3)}{, (#2#1#3)}{ an!d~(#2#1#3)}

\definecolor{3}{HTML}{FA0000}
\definecolor{7}{HTML}{008000}
\definecolor{13}{HTML}{0000A0}
\definecolor{17}{HTML}{FA00FA}
\definecolor{18}{HTML}{00F0FA}
\definecolor{14}{HTML}{903400}
\definecolor{19}{HTML}{F00000}
\definecolor{20}{HTML}{F00000}
\definecolor{8}{HTML}{4d8a9e}
\definecolor{9}{HTML}{4d8a9e}
\definecolor{10}{HTML}{e89530}
\definecolor{11}{HTML}{ffc885}
\definecolor{12}{HTML}{49006e}
\definecolor{5}{HTML}{7f9450}
\definecolor{4}{HTML}{ff8614}
\definecolor{6}{HTML}{b4eb34}
\definecolor{16}{HTML}{813feb}
\definecolor{15}{HTML}{00abc9}

\usepackage{tikz}
\usetikzlibrary{quotes,calc,positioning,arrows,arrows.meta,backgrounds,patterns}
\tikzset{
	cell complex/.style={
	filling/.style={
		fill=#1, fill opacity=.45
	},
	every node/.style={
		inner sep=1pt,
		fill=white,
		fill opacity=1,
		rounded corners=3pt
	},
	halo/.style={
		preaction={
			draw,
			white,
			line width=2pt,
			-
		}
	},
	every edge quotes/.append style={
		fill=none,
		inner sep=1,
		font=\scriptsize
	},
	t3/.style={filling=3},
	t4/.style={filling=4},
	t5/.style={filling=5},
	t6/.style={filling=6},
	t7/.style={filling=7},
	t8/.style={filling=8},
	t10/.style={filling=10},
	t12/.style={filling=12},
	t13/.style={semithick, draw=13, every edge quotes/.append style={text=13}},
	t14/.style={semithick, draw=14, every edge quotes/.append style={text=14}},
	t15/.style={semithick, draw=15, every edge quotes/.append style={text=15}},
	t16/.style={semithick, draw=16, every edge quotes/.append style={text=16}},
	t17/.style={preaction={draw, magenta!30, ultra thick, -Triangle Cap, shorten >=2pt}, t13},
	t18/.style={cyan, semithick},
	ll/.style={font=\sffamily\tiny},
}}

\usepackage[date=year]{biblatex}
\DeclareSourcemap{
	\maps[datatype=bibtex]{
		\map{
			\pertype{article}
			\pertype{inproceedings}
			\step[fieldsource=doi,final]
			\step[fieldset=url,null]
			\step[fieldset=issn,null]
			\step[fieldset=isbn,null]
		}
	}
}

\theoremstyle{definition}
\newtheorem*{definition*}{Definition}

\newcommand\F{\mathbf{F}}
\newcommand{\R}{\mathbf{R}}
\newcommand{\Z}{\mathbf{Z}}
\newcommand{\Cfg}{\mathit{Cfg}}
\newcommand{\SE}{\mathrm{SE}}
\newcommand{\id}{\mathrm{id}}
\mathchardef\mhyphen="2D
\newcommand{\Vect}[1]{#1\mhyphen\mathrm{Vec}}
\newcommand{\VR}{\mathit{VR}}
\newcommand{\RP}{\mathbf{P}}

\newlength\TempLenA 

\begin{document}
\maketitle
\begin{abstract}
	\noindent
	The conformation space of cyclooctane,
	a ringlike organic molecule comprising eight carbon atoms,
	is a two-dimensional algebraic variety,
	which has been studied extensively for more than 90 years.
	We propose a cell structure representing this space,
	which arises naturally by partitioning the space
	into subsets of conformations that admit particular symmetries.
	We do so both for the labeled conformation space,
	in which the carbon atoms are considered as distinct,
	and for the actual, unlabeled, conformation space.
	The proposed cell structure is obtained
	by identifying subspaces of conformations based on symmetry patterns
	and studying the geometry and topology of these subsets
	using methods from dimensionality reduction and topological data analysis.
	Our findings suggest that, in contrast to the labeled variant,
	the conformation space of cyclooctane is contractible.
\end{abstract}
\tikzset{cell complex}

\tableofcontents

\section{Introduction}
The fact that most alkanes are not rigid molecules but admit some flexibility
has caused a certain interest in understanding
the geometry of the space of all “shapes”, called \emph{conformations}, that a given alkane can attain.
Cycloalkanes are molecules forming a closed loop of carbon atoms,
which restricts their flexibility considerably
and gives their conformation spaces a particularly rich structure.
In particular, conformation spaces of small cycloalkanes (with five to around twelve carbon atoms) have been the subject
of various experimental and theoretical studies.
In the present work, we consider cyclooctane,
an alkane consisting of a cycle of eight carbon atoms,
each of which binds to two hydrogen atoms.

The conformations that can be attained by a flexible molecule
are governed by the potential energy required to form them,
which results mostly from
deviations of the bond lengths
and bond angles from their equilibria,
as well as steric and electronic interactions between hydrogen atoms bound to different carbon atoms.
The terms representing bond lengths and bond angles
grow by several orders of magnitude faster than the others;
therefore, it is customary to regard these quantities as fixed
at the equilibria \cite{Hendrickson:1967,GibsonScheraga:1997}.

Since in the case of cyclooctane the bond angles do not allow for the formation of a planar, regular octagon,
any conformation necessarily has to “pucker” in three dimensions,
meaning that the molecule has to fold and twist into a non-flat shape in order to accomodate the bond angles in three dimensional space;
see e.g.\ \cref{fig:symmetry-types} for a few typical conformations.
Because of this, conformations of this molecule are not rigid
but have a certain amount of mobility,
which gives rise to a rich structure of the conformation space.

\subsection{Representing molecules as linkages}
As a simplified model, we regard cyclooctane as an octagon embedded into three-dimensional Euclidean space
with all edges having the same fixed length
and all pairs of adjacent edges enclosing the same fixed angle.
Such an object is an instance of a \emph{linkage}:

\begin{definition*}
	An \emph{abstract linkage} is a pair $(K, \ell)$
	consisting of a graph $K$  with vertices $K_0$ and edges $K_1 \subseteq K_0 \times K_0$,
	and a length function $\ell\colon K_1\to (0,\infty), (i,j) \mapsto \ell_{ij}$ such that $\ell_{ij} = \ell_{ji}$ for all $i$, $j$.
	A \emph{realization} of $(K, \ell)$ in the Euclidean space $\mathbf R^n$
	is a mapping $x\colon K_0\to \mathbf R^n, i \mapsto x_i$
	such that for all edges $(i,j)\in K_1$ the Euclidean distance $\lVert x_i - x_j\rVert$ satisfies $d(x_i, x_j) = \ell_{ij}$.

	An \emph{abstract linkage with fixed joints} $(K, \ell, \phi)$
	is an abstract linkage endowed with a map $\phi\colon\{(i,j,k)\mid (i,j), (j,k) \in K_1\}\to \R/2\pi\Z, (i,j,k) \mapsto \phi_{ijk}$,
	such that $\phi_{ijk} = \phi_{kji}$ for all $i$, $j$, $k$.
	The map $\phi$ can be thought of as assigning an angle to each pair of adjacent edges.
	A \emph{realization} of $(K, \ell, \phi)$ is a realization $x$ of $(K, \ell)$
	such that $\angle(x_i, x_j, x_k)=\phi_{ijk}$ for all $(i,j), (j,k) \in K_1$.
		We write $R^n(K, \ell, \phi)$ or just $R^n(K)$ for the space of realizations in $\mathbf R^n$.
\end{definition*}

The study of linkages goes back at least to \cite{Kempe:1877};
their connection to algebraic sets has been studied in \cite{JordanSteiner:1999}.
Note that an abstract linkage with fixed joints can equivalently be described as an abstract linkage (without joints),
by additionally prescribing the distances $d(x_{i}, x_{k})$ for $(i,j), (j,k) \in K_1$.
Specifically, by the cosine law,
realizations of $(K, \ell, \phi)$ satisfy the quadratic equations
\begin{equation*}
	\begin{alignedat}{2}
		\lVert x_i - x_j\lVert^2 &= \ell_{ij}^2\quad && \text{for all}\ (i,j)\in K_1,\\
		\lVert x_i - x_k\lVert^2 &= \ell_{ij}^2 + \ell_{jk}^2 - 2\ell_{ij}\ell_{jk}\cos\phi_{ijk}\quad &&\text{for all}\ (i,j), (j,k)\in K_1.
	\end{alignedat}
\end{equation*}
Therefore, $R^n(K, \ell, \phi)$ is a real algebraic variety
(not necessarily connected, see \cite{Membrillo-SolisPirashviliEtAl:2019}).
The special Euclidean group $\SE(n)$,
generated by translations and rotations in $\mathbf R^n$,
acts on realizations by composition;
given a realization $x$ of $(K, \ell, \phi)$, any $T\in \SE(n)$
yields another realization $Tx\coloneqq T\circ x\colon K_0\to \mathbf R^n$ of $(K, \ell, \phi)$.

\begin{definition*}
	The quotient $\Cfg^n(K, \ell, \phi)\coloneqq R^n(K, \ell, \phi)/\SE(n)$ (or just $\Cfg^n(K)$)
	is called the \emph{configuration space} of $K$.
	Its elements are called \emph{configurations}.
	We write $\bar{x}$ for the equivalence class in $\Cfg^n(K)$ of a configuration $x \in R^n(K)$.
\end{definition*}

To see that $\Cfg^n(K, \ell, \phi)$ is also a variety,
choose three vertices $i, j, k\in K_0$
and fix three points $x_i, x_j, x_k\in\mathbf R^3$
that satisfy $\lVert x_i-x_j\rVert = \ell_{ij}$, $\lVert x_j-x_k\rVert = \ell_{jk}$ and $\angle(x_i, x_j, x_k)=\phi_{ijk}$.
Then every equivalence class in $\Cfg^n(K, \ell, \phi)$ has a unique representative $x$ with the prescribed images $x_i$, $x_j$ and $x_k$.

\subsection{Cyclooctane}
From now on, we consider the closed loop linkage $(K, \phi)$
with eight vertices $K_0=\{v_0,\dotsc,v_7\}$,
edges $e_0=(v_0,v_1),\dotsc,e_6=(v_6,v_7),=e_7(v_7,v_0)$ of equal length $\ell = \ell_{i,i+1}$
and fixed joint angles $\phi=\phi_{i-1,i,i+1}$ for all $i \in K_0$, where we consider all indices modulo~8.

\begin{wrapfigure}{O}{3cm}
	\centering
	\includegraphics[width=3cm]{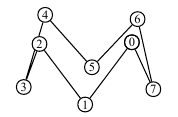}%
	\caption{A typical configucation of the cyclooctane linkage}
	\label{fig:5223}
\end{wrapfigure}
We are only interested in configurations in three dimensions,
and write $R(K)$ for $R^3(K)$ and $\Cfg(K)$ for $\Cfg^3(K)$.
For $\phi>\frac{3\pi}{4}$, $\Cfg(K)$ is empty,
and for $\phi=\frac{3\pi}{4}$, $\Cfg(K)$ has a single element,
the embedding of a planar regular octagon into $\mathbf R^3$.
For $\phi<\frac{3\pi}{4}$,
a configuration necessarily folds and twists in all three spatial dimensions; see \cref{fig:5223} and, for more examples, \cref{fig:symmetry-types}.
An actual cyclooctane molecule has
$\ell \approx 1.53\,\mathrm{A}$ and $\phi\approx 115^\circ$
\cites[Table~II]{MeiboomHewittEtAl:1977}[Table~I]{Hendrickson:1967}.
This is the angle $\phi$ in the sample of $\Cfg(K)$ used in our experiments \cite{Martin:2012}.

The space $\Cfg(K)$ has attracted a large attention in the literature
because it is a two-dimensional variety.
Hence, it is large enough to have interesting structure and topology,
while still being manageably low-dimensional for manual scrutiny.

\subsection{Labeled and unlabeled configurations}
\label{sec:labeled-and-unlabeled}
The map $x$ keeps track of the individual images of the vertices in $K$;
the space $\Cfg(K)$ is in fact the space of all \emph{labeled} configurations.
However, as all carbon atoms in an actual molecule are indistinguishable,
the labeling is only a technical artifact, and
we are actually primarily interested in the topological structure of the quotient of $\Cfg(K)$
by the symmetry group of the regular octagon
(i.e., the dihedral group $D_8=C_8\rtimes C_2$, where $C_8$ acts by rotation and $C_2$ by reversal of the indexation),
together with the quotient map $\Cfg(K) \to \Cfg(K)/D_8$.
As an intermediate step, we also analyze the topology of the quotient of $\Cfg(K)$
by the subgroup $C_8$ of $D_8$,
together with the quotient maps $\Cfg(K) \to \Cfg(K)/C_8$ and $\Cfg(K)/C_8 \to \Cfg(K)/D_8$.

We note that, surprisingly, the unlabeled conformation space of cyclooctane seems not to have been discussed in the literature to the best of our knowledge, even though the study of equivalence classes of alkanes up to relabeling of carbon atoms is a classical topic dating back to Pólya \cite{Polya:1937}.

\subsection{Related work}
The conformation space of cyclooctane has been subject to a large number of studies.
We discuss the literature most relevant to our work.

\runinhead{Conformational energies}
From the perspective of molecular chemistry,
it is interesting to understand not only the space of realizable conformations,
but also the function $E\colon \Cfg\to[0,\infty)$
that assigns to each conformation its conformational energy.

In his influential work \cite{Hendrickson:1967},
\citeauthor{Hendrickson:1967} discusses various terms
contributing to the total conformational energy.
With regard to their high stiffness, he assumes bond lengths to be fixed for all his calculations;
see also \cite{GibsonScheraga:1997}.
In \cite{Hendrickson:1967b},
the shapes of conformations of locally minimal conformational energy (there called \emph{conformers})
are described in terms of their symmetries as well as their dihedral angles,
and the names under which they are referred to.
\Citeauthor{Hendrickson:1967b} also provides a nomenclature for certain symmetries in the conformations.
He describes the symmetric, asymmetric and pseudo-rotational modes
of interconversion between the conformations.

Since then, research has mostly focused on these conformers.
Various different approaches for the description of conformers have been proposed,
leading to the development of so-called \emph{ring puckering} coordinates \cite{CremerPople:1975};
these describe a conformer by the out-of-plane displacement
in terms of an amplitude and a phase coordinate  w.r.t\ canonical conformations.
Another approach \cite{PickettStrauss:1971} in terms of irreducible representations
of dihedral groups (the symmetry groups of regular polygons), equivalent to ring puckering coordinates \cite{BoeyensEvans:1989}, have proven fruitful for conformational analysis \cite{BoeyensEvans:1989,EvansBoeyens:1989}.

\runinhead{Sampling}
Sampling algebraic varieties such as $\Cfg(K)$ is an involved problem in its own right.
A feasible method should yield a uniform, sufficiently dense sample (i.e., solutions from the entire solution space are found)
and allow for the estimation of the uncertainty of the solutions if necessary.

\runinhead{Cayley–Menger determinants}
A sufficient and necessary condition
for the embedability of finitely many points into a Euclidean space with prescribed distances
is given by the signs of the \emph{Cayley--Menger-determinants}
\cite[corollary 42.2]{Blumenthal:1970}.
A technique based on Cayley--Menger determinants is used in \cite{ThomasPortaEtAl:2004,PortaRosEtAl:2007}
to develop a sampling technique that is used for sampling the configuration spaces
describing conformations of cyclooctane and other short cyclic molecules.

\runinhead{Resultants}
\Citeauthor{CoutsiasSeokEtAl:2006} develop a sampling algorithm using Sylvester and Dixon resultants
\cite{CoutsiasSeokEtAl:2006}.
These are successfully employed in \cite{MartinThompsonEtAl:2010}
to compute a sample of the cyclooctane conformation space
of approximately one million configurations.
This sample, available online \cite{Martin:2012},
also includes positions of the hydrogen atoms for the actual cyclooctane molecule;
as there are two hydrogen atoms per carbon atom, the full conformations are
thus given as a collection of points in $\mathbf R^{72}$
(eight carbon atoms sixteen hydrogen atoms with three spatial coordinates each).
Note however that the positions of the hydrogen atoms are uniquely determined by those
of the carbon atoms.
All conformations in this sample have their center of mass at the origin
and are aligned according to the Eckart condition (see \cref{sec:defining-metrics}).
This is the sample employed in our analysis.

\runinhead{Homotopy continuation}
Recently, homotopy continuation has proven to be a worthwhile method
to tackle computational algebro-geometrical problems
\cite{SommeseWamplerCharlesWeldon:2005}.
Implementations for various environments are available;
the Julia package by \citeauthor{BreidingTimme:2018} \cite{BreidingTimme:2018}
provides code for uniformly sampling the cyclooctane conformation space.

\runinhead{Low-dimensional projections, surface reconstruction, topology}
Among the more recent work on cyclooctane,
the investigations carried out in \cite{BrownMartinEtAl:2008,MartinThompsonEtAl:2010,MartinWatson:2011,Martin:2012}
are of particular interest for us.
Based on the incremental surface reconstruction algorithm in \cite{Freedman:2007},
which infers a triangulation of an orientable or non-orientable 2-manifold
from a finite number of points,
\citeauthor{MartinWatson:2011}
presented an algorithm \cite{MartinWatson:2011}
that allows for the reconstruction of a 2-dimensional stratified space with singularities.
The feasibility of their method is demonstrated both in terms of running times
as well as quality of the computed triangulations,
although neither \cite{MartinWatson:2011} nor \cite{Freedman:2007}
provide complexity bounds or qualitative guarantees.

In \cite{MartinThompsonEtAl:2010},
the method is applied to the cyclooctane sample from \cite{Martin:2012}.
The result is compared to the image of the conformation space
under a projection to three dimensions obtained using \emph{Isomap},
a non-linear dimensionality reduction method (see \cref{sec:Isomap} and \cite{Tenenbaum:2000}).
In addition, the reconstruction of the actual surface enables the authors to
identify the conformation space as a gluing of a sphere and a Klein bottle.
As far as we know, \cite{MartinThompsonEtAl:2010} is the earliest work
suggesting this description of the conformation space of cyclooctane.
Their cyclooctane data set has also been studied using persistent homology in \cite{AdamsMoy:2021}.
We extend the observations made there.

\runinhead{Persistent homology and Morse theory}
Persistent homology \cite{EdelsbrunnerHarer:2008,ZomorodianCarlsson:2005,Ghrist:2007} is a tool from algebraic topology that describes topological features of a filtration.
It is often applied to geometric filtrations such as the Vietoris--Rips filtration, which can be used to infer
(under mild assumptions) the actual homology of a sampled topological space.
\Citeauthor{Membrillo-SolisPirashviliEtAl:2019} undertake a two-step analysis of the conformation spaces of cyclooctane and other molecules,
and the conformational energy distribution on these.
In their analysis \cite{Membrillo-SolisPirashviliEtAl:2019},
which is also performed on the sample from \cite{Martin:2012},
they use persistent homology to validate hypotheses about the
homology of several conformation spaces.
They further employ a numerical approximation of the Morse--Smale complex
\cite{GyulassyBremerEtAl:2011,EdelsbrunnerHarerEtAl:2003,GyulassyBremerEtAl:2008} of the conformational energy function based on a triangulation of the data set,
applying topological simplification guided by persistent homology \cite{GuntherReininghausEtAl:2014,BauerLangeEtAl:2012} to extract the significant topological features of the energy function.

\runinhead{Symmetries}
An early study of conformation spaces of molecules
in terms of actions of symmetry groups is carried out in \cite{BoeyensEvans:1989,PickettStrauss:1971}.
A thorough investigation of idealized configuration spaces of linkages (there called molecular graphs)
and their quotients w.r.t.\ their symmetry groups is done in \cite{Membrillo-SolisPirashviliEtAl:2019}.
The authors observe that the symmetry group actions
render the respective quotients of the configuration spaces into orbifolds.
The work of \cite{Membrillo-SolisPirashviliEtAl:2019} appears to be the first
to introduce other metrics than the Euclidean metric on $\Cfg(K)$,
which are defined in terms of orbifold theory.

\subsection{Our contribution}
We provide further insight into the structure of the conformation space of cyclooctane, supporting the description
as the gluing of a sphere and a Klein bottle.
Our findings are based on a partition of $\Cfg(K)$ into certain subspaces that give rise to a cell structure on $\Cfg(K)$.
Each of these subspaces arises naturally as the set of all configurations in $\Cfg(K)$ that exhibit a common \emph{symmetry type}; see \cref{sec:symmetry-types}.

Symmetry types of cycloalkanes \cite{Hendrickson:1967,Hendrickson:1967b}
have already been seen to distinguish conformations from the sphere and the Klein bottle subspace of $\Cfg(K)$ \cite{MartinThompsonEtAl:2010}.
We extend this approach significantly, identifying more symmetry types
and using their regimes for the segmentation of $\Cfg(K)$ into a cell complex.
As far as we know, no cell complex representing the conformation space of cyclooctane has been proposed so far.

Besides the labeled configuration space $\Cfg(K)$, we also provide a description of the unlabeled configuration space $\Cfg(K)/D_8$,
which arises as a quotient of the labeled configuration space under the action of the symmetry group $D_8$.
Interestingly, the latter space has a much simpler topological structure;
namely, while $\Cfg(K)$ is obtained as the gluing of a sphere and a Klein bottle,
the space $\Cfg(K)/D_8$ is contractible.

The experiments and the resulting plots have been made
using a 6040-point subsample of the data set \cite{Martin:2012}%
\footnote{%
	See the download “cyclo-octane example” in section “NMTRI: Non-Manifold Surface Reconstruction” on their website.
	For details, consider the file \texttt{readme.txt} in the download.
}.
The subsample is provided as part of the data set, together with code to generate it, and subsampled to ensure a minimum
distance between data points.

\section{Methods}
\subsection{Metrics}
\label{sec:defining-metrics}
We consider two metrics on $\Cfg(K)$.

\runinhead{Angular metric}
Let $\bar{x} \in \Cfg(K)$ be a configuration, represented by the realization $x \in R(K)$.
We write $e_i$ for the vector $e_i \coloneqq x_{i+1} - x_i$.
To $x$, we assign a sequence $\sigma(x)=(\sigma(x)_0,\dotsc,\sigma(x)_7) \subset (S^1)^8$
of oriented angles, called the \emph{torsion angles} or \emph{dihedral angles} of $x$, as follows.

Using $\|\cdot\|$ to denote the Euclidean norm
and writing
$n_{i} = \frac{e_i \times e_{i+1}}{\lVert e_i \times e_{i+1}\rVert}$
(with indices taken modulo $8$),
let $\sigma_i \coloneqq \angle(n_{i-1}, n_{i})$ be the oriented angle such that $\sigma_i$ is positive if $(n_{i-1}, n_i, e_i)$ is a right handed system;
that is, if $\langle n_{i-1} \times n_i, e_i\rangle$ is positive.
Therefore,
\[
	\sigma_i = \operatorname{atan2}\bigl( \langle n_{i-1} \times n_i, e_i\rangle, \langle n_{i-1}, n_i\rangle \bigr),
\]
where $\operatorname{atan2}$ is the continuous extension of $(a,b)\mapsto\arctan\frac ab$.
Torsion angles are invariant under the action of $\SE(3)$,
so the definition of $\sigma(x)$ passes down to a well-defined map $\sigma\colon \Cfg(K) \to (S^1)^8$, $\bar{x} \mapsto \sigma(x)$.
We call the sequence $\sigma(\bar{x})$ the \emph{torsion} or \emph{dihedral angles} of the configuration $\bar{x}$.
Moreover, two configurations agree if and only if all their torsion angles do.
Therefore, letting
\[
	d_\angle(\bar{x}, \bar{x}') = \sum_i \min\bigl(\lvert \sigma(x)_i - \sigma(x')_i \rvert, 2\pi - \lvert \sigma(x)_i - \sigma(x')_i \rvert\bigr)
\]
defines a metric on $\Cfg(K)$,
called the \emph{angular metric}.

\runinhead{Euclidean metric}
\begin{wrapfigure}{o}{3cm}
	\centering
	\includegraphics[width=3cm]{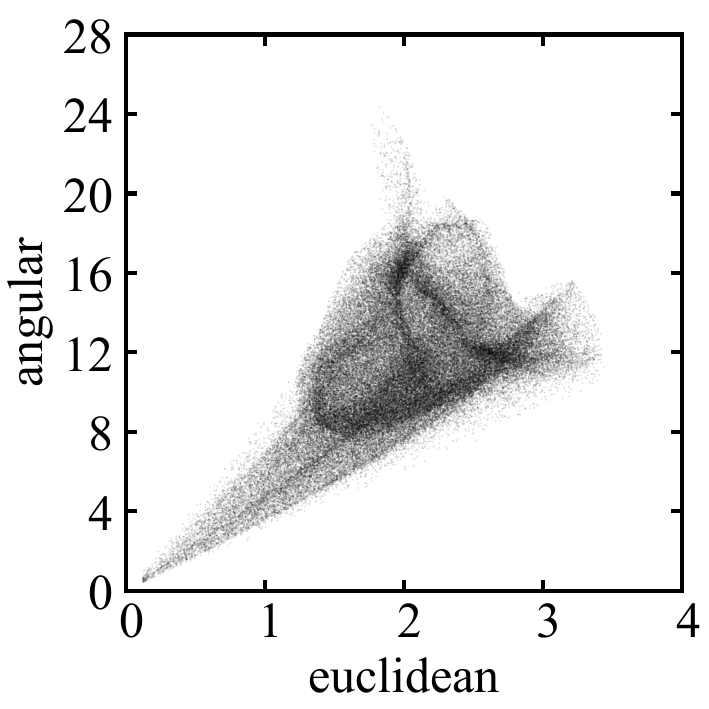} 
	\caption{%
		The Euclidean distance $d_\Vert$ and the angular distance $d_\angle$
		are equivalent.
		The picture has been obtained by plotting 100k randomly choosen pairwise distances (out of the $\approx 18M$ possible ones).
	}
	\label{fig:distance_correlation}
\end{wrapfigure}
In the sample \cite{Martin:2012} that we used in our analysis,
each equivalence class in $\Cfg(K)$ is represented by a \emph{standard representative} that is aligned as follows.
Let $x^\mathrm{pln}_i$ be the coordinates of the planar regular octagon in the $xy$-plane with circumradius~1,
given by
\( x^\mathrm{pln}_i = (\cos (2i-1)\frac{\pi}{8}, \sin (2i-1)\frac{\pi}{8}, 0 )^\mathrm T \).
For the equivalence class $\bar{x}$ of a realization $x\in R(K)$ of $K$,
\citeauthor{Martin:2012} chooses the representative $x^\mathrm{std}$ of $\bar{x}$ that satisfies the \emph{Eckart condition}
\begin{equation}
	\label{eq:Eckart condition} \sum_i x^\mathrm{std}_i = 0, \qquad \sum_i x^\mathrm{pln}_i \times x^\mathrm{std}_i = 0
\end{equation}
w.r.t.\ $x^\mathrm{pln}$ \cite{Eckart:1935}.
The \emph{Euclidean metric} is the metric $d_{\Vert}$ on $\Cfg(K)$
defined as
\[d_{\Vert }(x, x')^2=\sum_i \lVert x^{\text{std}}_i - x'^{\text{std}}_i\rVert^2\]
for two configurations $\bar{x}$ and $\bar{x}'$ in $\Cfg(K)$,
where $\lVert \cdot\rVert$ denotes the usual Euclidean metric on $\R^3$.

The angle sequence $(\sigma_0,\dotsc,\sigma_7)$ uniquely determines a equivalence class in $\Cfg(K)$.
A method for sampling $\Cfg(K)$ in terms of squared distances $\Vert x_{i-1}-x_{i+2}\Vert^2$,
and thus in terms of dihedral angles, is described in \cite{PortaRosEtAl:2007}.

Since $d_\angle$ is invariant under the action of $\SE(3)$ and involves no alignment to a reference orientation,
one could consider it a more “natural” metric than $d_{\Vert}$.
However, since previous investigations have considered $d_{\Vert}$
as the only choice of a metric on $\Cfg(K)$,
we include it in our investigations nevertheless.
As it turns out, both metrics are equivalent with low metric distortion, as can be seen from \cref{fig:distance_correlation}.

\subsection{Persistent homology}
\label{sec:persistent-homology}
We use \emph{persistent homology} \cite{EdelsbrunnerHarer:2008} as a central tool in our analysis.
In the following, we give a quick account of persistent homology; for details, see \cite{EdelsbrunnerHarer:2008,OtterPorterEtAl:2017}.
Given a finite metric space $(X=\{x^1,\dotsc, x^k\}, d)$,
the associated \emph{Vietoris–Rips complex} $\VR_*(X)$
assigns to each $r\in[0, \infty)$
the simplicial complex $\VR_r(X)$ with $\VR_r(X) = \{\sigma \subseteq X \mid \emptyset \neq \sigma, \ \forall x, x' \in \sigma\colon d(x, x') \leq r \}$,
and to every pair $r \leq s$ in $[0, \infty)$
the inclusion map $\rho_{rs}\colon\mathit{VR}_r(X) \hookrightarrow \mathit{VR}_s(X)$.

A persistence module is a functor $[0, \infty) \to \Vect{\F}$,
to the category of vector spaces over a fixed field $\F$.
Explicitly, a persistence module $M$ assigns to every $r \in [0, \infty)$ an $\F$-vector space $M$
and to every pair $r \leq s$ in $[0, \infty)$ a map $M_{sr}\colon M_r \to M_s$,
such that $M_{rr} = \id$ and $M_{sr} M_{ts} = M_{tr}$ for all $r \leq s \leq t$ in $[0, \infty)$.

If $M$ is a persistence module such that $\dim M_r < \infty$ for all $r$,
then there is an essentially unique decomposition $M \cong \bigoplus_{i \in I} \F^{[b_i, d_i)}_*$
of $M$ into a direct sum of \emph{interval modules} 
\cite{Crawley-Boevey:2015}, defined as follows.
For $0\leq b<d\leq\infty$,
the \emph{interval module} $\F^{[b,d)}$ supported on $[b, d)$ is the persistence module with
\[
	\F^{[b,d)}_r = \begin{cases}
		\F&\text{if $b \leq r < d$}\\
		0&\text{otherwise,}
	\end{cases}
	\qquad
	\F^{[b,d)}_{sr} = \begin{cases}
		\id&\text{if $b \leq r \leq s < d$}\\
		0&\text{otherwise}
	\end{cases}
\]
for all $r$, $s$.
Taking the $q$-th simplicial homology $H_q(X, \F) \coloneqq H_q(\VR_*(X), \F)$ of $\VR_*(X)$ with coefficients in $\F$ defines a persistence module,
called the $q$-th \emph{persistent homology} of $\VR_*(X)$.
Informally, the intervals $[b_i, d_i)$ in the decomposition $H_q(X, \F) \cong \bigoplus_{i \in I} \F^{[b_i, d_i)}_*$ of $H_q(X, \F)$ capture the ranges of parameters $r$
for which a specific homological feature persists.
Note that $H_q(X, \F)$ depends on the field $\F$.

The intervals $[b_i, d_i)$ occurring in the decomposition $H_q(X, \F) \cong \bigoplus_{i \in I} \F^{[b_i, d_i)}_*$
are commonly drawn as dots in a \emph{persistence diagram}, where the axes correspond to the endpoints of the intervals.
In practice, one cuts off the Vietoris--Rips complex at some threshold $r_\infty$.
Points with $d_i \geq r_\infty$ are drawn on the upper boundary.
\begin{figure}[tbp]
	\centering
	\subcaptionbox{$d_\angle$}{\includegraphics{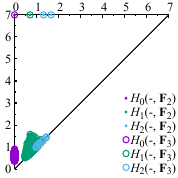}}
	\subcaptionbox{$d_\Vert$}{\includegraphics{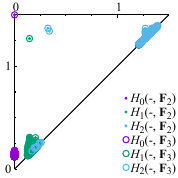}}
	\caption{%
		Persistent homology of $\Cfg(K)$ with coefficients in $\mathbf{F}_2$ and $\mathbf{F}_3$ w.r.t.\ $d_\angle$ and $d_\Vert$.
		Each point $(b, d)$ corresponds to one homology interval $[b, d)$.
		For this sample, the diagrams for $H_\bullet(-, \F_2)$ and $H_\bullet(-, \F_3)$ coincide;
		this is not the case for all subspaces we consider.
	}
	\label{fig:persistence:config-space}
\end{figure}
For instance, \cref{fig:persistence:config-space} shows the persistent homology
of $\Cfg(K)$ with respect to the two metrics $d_\angle$ and $d_\Vert$ with coefficients in the finite fields $\F_2$ and $\F_3$.

The applicability of persistent homology is justified by the following statement:
If $M$ is a manifold and $X\subseteq M$ is a sufficiently large i.i.d.\ sample from $M$,
then the long living intervals in $H_\bullet(\VR_*(X))$
capture the actual singular homology $H_\bullet(M)$ of $M$ with high probability
\cites[Thm.\ 3.1]{NiyogiSmaleEtAl:2008}[Thm.\ 10]{AttaliLieutierEtAl:2013}.
See \cite{EdelsbrunnerHarer:2008,EdelsbrunnerMorozov:2017,OtterPorterEtAl:2017} for more details on persistent homology.

Several implementations are available for computing persistent homology of Vietoris–Rips complexes
from a given distance matrix (see \cite{OtterPorterEtAl:2017} for an overview).
Because of its computational efficiency,
we use the freely available software Ripser \cite{Bauer:2021}
which can compute the interval endpoints with respect to coefficients in finite fields.

For example, we learn the following about $\Cfg(K)$ from \cref{fig:persistence:config-space}:
In the bottom left corner of the diagrams, near the diagonal, we see many short-lived topological features,
while the four points in the top left represent long-lived features.
In fact, the diagram shows that there is a long range of filtration parameters (from 2 to at least 7 for $d_\angle$)
during which the homology type of the Vietoris--Rips complex does not change.
One may thus argue that with high probability, the four points in the top left represent the true homology of $\Cfg(K)$,
so that that $\Cfg(K)$ has the $\F_2$- and $\F_3$-Betti numbers $(1,1,2)$ in dimensions $(0,1,2)$, respectively.

Computation of persistent homology depends on the field $\F$.
We draw all our diagrams both for $\F_2$ and $\F_3$-valued persistent homology.
Although both coincide in \cref{fig:persistence:config-space}, we will see subspaces of $\Cfg(K)$ where they do not.
\phantomsection\label{sec:classification of surfaces}
We do so because a closed surface (a two-dimensional compact manifold without boundary) is characterized uniquely up to homeomorphism by its homology with coefficients in $\F_2$ and $\F_3$, since closed surfaces are classified by their orientability and their genus, and homology over $\F_3$ detects orientability, while homology over $\F_2$ identifies the (orientable or non-orientable) genus.

\subsection{Isomap}
\label{sec:Isomap}
Using the coordinates of the vertices, the space $\Cfg(K)$ can be embedded into $\mathbf R^{24}$
(eight vertices, with three spatial coordinates each).
Analogously, using the torsion angles, it can be embedded into the $8$-fold product $S^1\times\cdots\times S^1$.

To obtain projections onto lower-dimensional Euclidean spaces, we apply a method called \emph{Isomap} \cite{Tenenbaum:2000}.
Given a finite subset $X = \{x_1,\dotsc,x_k\} \subseteq \R^N$ of supposedly low intrinsic dimension embedded in a high-dimensional Euclidean space,
this method first estimates an \emph{intrinsic metric} $d_X$ from the Euclidean metric $d_{\R^N}$ and then applies \emph{multidimensional scaling}
to compute a (non-linear) projection $\tilde{X} = \{\tilde{x}_1,\dotsc,\tilde{x}_k\} \subset \R^n$ of $X$ embedded in a low dimensional Euclidean space with $n \leq N$,
such that the restriction of the Euclidean metric $d_{\R^n}$ to $\tilde{X}$ approximates $d_X$.
With our Isomap projections, we also calculate the \emph{metric distortion} $\log \frac{d_{\R^n}(\tilde{x}_i, \tilde{x}_j)}{d_X(x_i, x_j)}$ of each pair $(x_i, x_j) \in G(X)$.
This quantity is large for points that are placed close together in $\tilde{X}$ that are actually distant in $X$.

\begin{figure}[tbp]
	\centering
	\subcaptionbox*{}{\includegraphics[width=2cm]{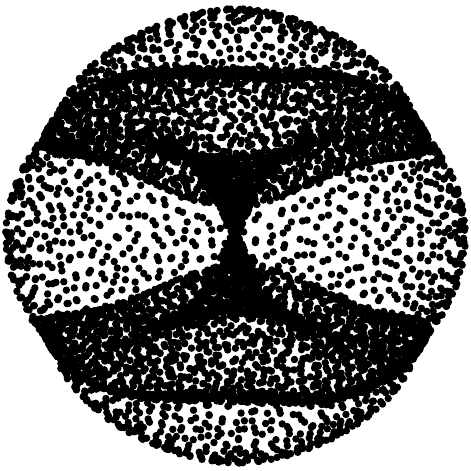}}%
	\raisebox{1cm-1ex}{${}   ={}$}\subcaptionbox*{}{\includegraphics[width=2cm]{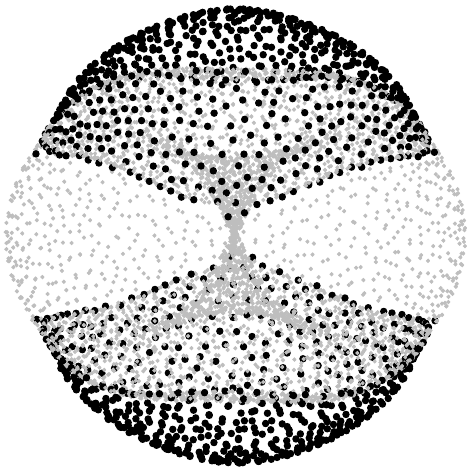}}%
	\raisebox{1cm-1ex}{${}\cup{}$}\subcaptionbox*{}{\includegraphics[width=2cm]{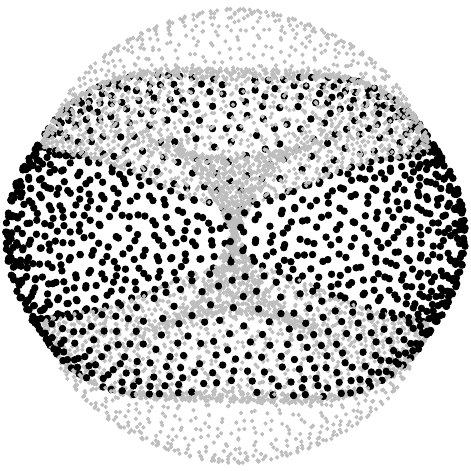}}%
	\raisebox{1cm-1ex}{${}\cup{}$}\subcaptionbox*{}{\includegraphics[width=2cm]{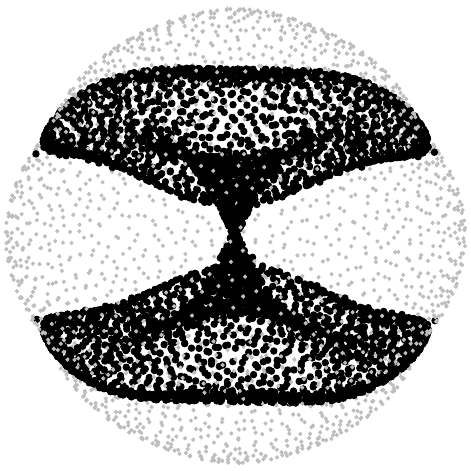}}\\[-\baselineskip]
	\subcaptionbox*{$\Cfg(K)$}{\includegraphics[width=2cm]{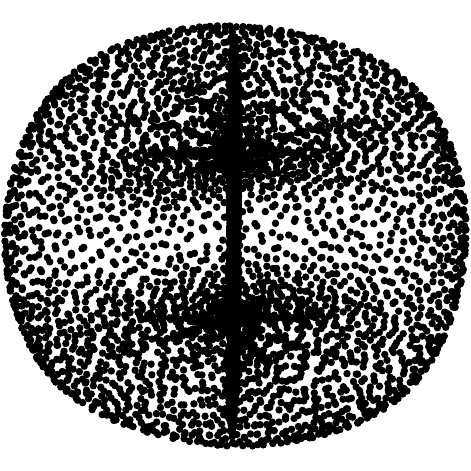}}%
	\raisebox{1cm-1ex}{${}   ={}$}\subcaptionbox*{$A$}{\includegraphics[width=2cm]{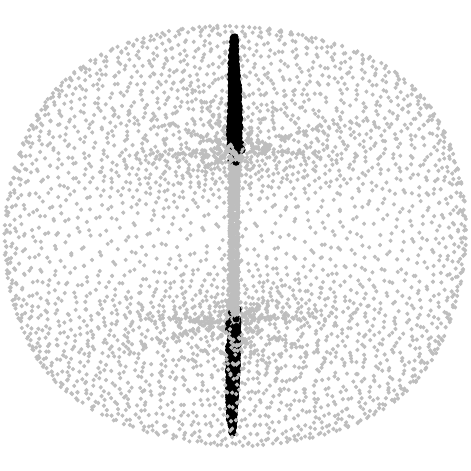}}%
	\raisebox{1cm-1ex}{${}\cup{}$}\subcaptionbox*{$B$}{\includegraphics[width=2cm]{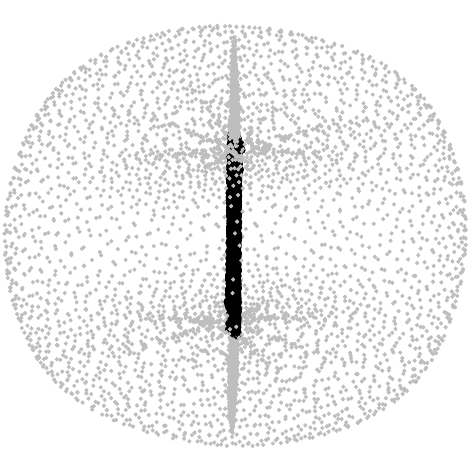}}%
	\raisebox{1cm-1ex}{${}\cup{}$}\subcaptionbox*{$C$}{\includegraphics[width=2cm]{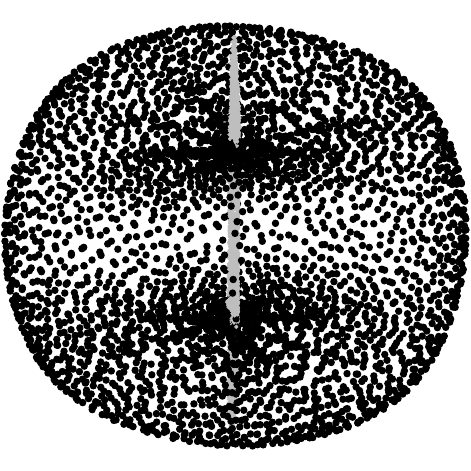}}\\
	\caption{Isomap projection (top: w.r.t.\ $d_\Vert$, bottom: $d_\angle$) of $\Cfg(K)$ and its subspaces $A$, $B$, $C$.}
	\label{fig:isomap subspaces}
\end{figure}

Isomap has already been applied to the cyclooctane data set in \cite{MartinThompsonEtAl:2010}.
In our case, we use the distances $d_\angle$ and $d_\Vert$ as input for Isomap
and obtain projections to $\R^3$ of $\Cfg(K)$ and various of its subsets.
An Isomap projection of $\Cfg(K)$ w.r.t.\ $d_\angle$ and $d_\Vert$ is shown in \cref{fig:isomap subspaces}.
We observe that the projections exhibit regions of high metric distortion around the two circles in which the outer sphere and the inscribed hourglass intersect.

\section{A cell structure for the labeled configuration space}
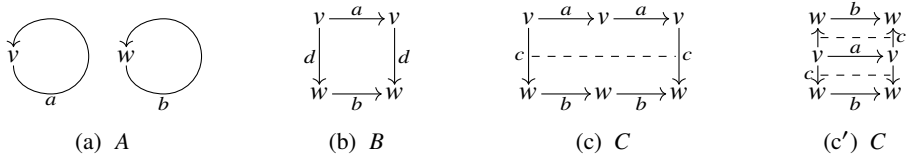
\begin{figure}[tbp]
	\tikzset{every edge/.append style={->}, every node/.append style={execute at begin node=$, execute at end node=$}}
	\subcaptionbox{$A$}{%
		\begin{tikzpicture}[baseline=(v)]
			\draw[shift={(0,0.5)}, ->]
				(0,0) node (v) {v}
				(0.5, -0.5) coordinate (t1)
				(1.0,  0.0) coordinate (t2)
				(0.5,  0.5) coordinate (t3)
				(v) to[out=south, in=west] (t1)
				to[out=east, in=south, "a"' at start] (t2)
				to[out=north, in=east] (t3)
				to[out=west, in=north] (v);
			\draw[shift={(1.5,0.5)}, ->]
				(0,0) node (w) {w}
				(0.5, -0.5) coordinate (t1)
				(1.0,  0.0) coordinate (t2)
				(0.5,  0.5) coordinate (t3)
				(w) to[out=south, in=west] (t1)
				to[out=east, in=south, "b"' at start] (t2)
				to[out=north, in=east] (t3)
				to[out=west, in=north] (w);
		\end{tikzpicture}
	}\hfill
	\subcaptionbox{$B$}{%
		\begin{tikzpicture}[baseline=(d)]
			\draw   (0,1) node (v1) {v}
				(1,1) node (v2) {v}
				(0,0) node (w1) {w}
				(1,0) node (w2) {w}
				(v1) edge["a" ] (v2)
				(v1) edge["d"'] (w1)
				(v2) edge["d" name=d] (w2)
				(w1) edge["b"'] (w2);
			\end{tikzpicture}
	}\hfill
	\subcaptionbox{$C$\label{fig:model:C}}{%
		\begin{tikzpicture}[baseline=(c)]
			\draw   (0,1) node (v1) {v}
				(1,1) node (v2) {v}
				(2,1) node (v3) {v}
				(0,0) node (w1) {w}
				(1,0) node (w2) {w}
				(2,0) node (w3) {w}
				(0,.5) node (t1) {}
				(2,.5) node (t3) {}
				(v1) edge["a" ] (v2) (v2) edge["a" ] (v3)
				(w1) edge["b"'] (w2) (w2) edge["b"'] (w3)
				(v1) edge["c"'] (w1) (v3) edge["c" name=c] (w3)
				(t1) edge[dashed, -] (t3);
		\end{tikzpicture}
	}\hfill
	\renewcommand{\thesubfigure}{\subref{fig:model:C}$'$}
	\subcaptionbox{$C$\label{fig:model:C:2}}{%
		\begin{tikzpicture}[y={(0,0.5)}, baseline={(v1)}]
			\draw (0,0)  node (v1) {v}
				(1,0)  node (v2) {v}
				(0,1)  node (w1) {w}
				(1,1)  node (w2) {w}
				(0,-1) node (w3) {w}
				(1,-1) node (w4) {w}
				(v1) edge["a"] (v2) edge coordinate (t1) (w1) edge["c"'] coordinate(t3) (w3)
				(v2) edge["c"'] coordinate (t2) (w2) edge coordinate(t4) (w4)
				(w1) edge["b"] (w2)
				(w3) edge["b"'] (w4)
				(t1) edge[dashed, -, shorten >=1pt, shorten <=1pt] (t2)
				(t3) edge[dashed, -, shorten >=1pt, shorten <=1pt] (t4);
		\end{tikzpicture}
	}
	\caption{Cellular model of the labeled configuration space $\Cfg(K)$, and its subspaces $A$, $B$ and $C$.
		The subspace $A$ consists of two disjoint disks,
		$B$ is a cylinder, and $C$ is a model for a Klein bottle sewed together by two Möbius strips
		(both with the dashed line as boundary; one with core curve $a$, one with core curve $b$; see (\subref{fig:model:C:2}) for a different drawing of the same model).
		The three subspaces are glued together along the 1-cells $a$ and $b$.
		Thus, the whole space is the union of the sphere $S^2 \cong A \cup B$ the Klein bottle $C$,
		where $A \cup B$ and $C$ intersect in the two disjoint circles $a$ and $b$.
	}
	\label{fig:cell complex}
\end{figure}

We will argue in the following that the labeled configuration space $\Cfg(K)$ this space is the union of three subspaces, denoted by $A$, $B$ and $C$,
such that $A$ consists of two disjoint 2-disks, $B$ is a cylinder (without the caps) and $C$ is a Klein bottle.
Thus, $\Cfg(K)$ is the union of the sphere $A \cup B$ and a Klein bottle $C$ that intersect in two disjoint circles (the 1-dimensional cells $a$ and $b$ in \cref{fig:cell complex}).
These are glued together as in \cref{fig:cell complex}.

This description has been previously provided in \cite{MartinThompsonEtAl:2010,MartinWatson:2011}.
We will provide additional evidence for this description.
More importantly, we will show that this cell complex arises naturally from the data.
Namely, we will see that the symmetry types (\cref{sec:symmetry-types}) endow $\Cfg(K)$ with a cell structure refining \cref{fig:cell complex}.

As a first justification of the proposed model, we see from the persistence diagram (\cref{fig:persistence:config-space})
that $\Cfg(K)$ has the Betti numbers $\beta_\bullet(\Cfg(K), \F_2) = \beta_\bullet(\Cfg(K), \F_3) = (1,1,2)$ in dimensions $(0,1,2)$, respectively,
which coincides with the Betti numbers of the proposed cell complex.

The Isomap projections of $\Cfg(K)$ with respect to the two metrics $d_\Vert$ and $d_\angle$ (\cref{fig:isomap subspaces}, left) serve as a starting point for our analysis.
Note that these projections are not isometric embeddings,
and do not necessarily reflect the topology of $\Cfg(K)$.

\subsection{Symmetry types}
\label{sec:symmetry-types}
We have not explained yet how the three subspaces $A$, $B$ and $C$ that partition $\Cfg(K)$ are defined.
A key ingredient to this end is the observation that the sequences of dihedral angles of the configurations in $\Cfg(K)$ exhibit certain \emph{symmetry types}.
In the following, we define certain symmetry types,
whose definitions have been chosen after scrutiny of the data,
that we will use in our analysis, and give numbers to them for easy reference; see \cref{tab:symmetry types}.
\begin{table}[tbp]
	{%
		\scriptsize
		\raisebox{.5\depth-0.5\height}{%
			\subcaptionbox{Two-dimensional strata, and the subspace $A$, $B$ or $C$ they belong to.\label{tab:2-dim types}}{
				\begin{tabular}[b]{c|c|rc}
					\toprule
					type & subspace & & condition \\ \midrule
					1 & $A$, $B$ & & $(\alpha, \beta, \gamma, \delta, \alpha, \beta, \gamma, \delta)$\\
					2 & $C$ & & not $(\alpha, \beta, \gamma, \delta, \alpha, \beta, \gamma, \delta)$\\\midrule
					3 & \multirow{2}{*}{$A$} &    & $\begin{psmallmatrix} \alpha, & \beta, & \gamma, & \delta, & \alpha, & \beta, & \gamma, & \delta \\ + & - & + & - & + & - & + & - \end{psmallmatrix}$\\
					4 & &    & $\begin{psmallmatrix}
						\alpha, & \beta & {}\lessgtr{} & \gamma & {}\gtrless{} & \delta, & \alpha, & \beta & {}\lessgtr{} & \gamma & {}\gtrless{} & \delta \\
						\pm     & \mp   &              & \mp    &              & \mp     & \pm     & \mp   &              & \mp    &              & \mp
					\end{psmallmatrix}$\\\midrule
					\multirow{2}{*}{5} & \multirow{4}{*}{$B$} &   & $\begin{psmallmatrix} \alpha, & \beta &{}\lessgtr{} & \gamma &{}\lessgtr{}&\delta,& \alpha, & \beta &{}\lessgtr{} & \gamma &{}\lessgtr{}&\delta \\ \pm & \mp && \mp && \mp & \pm & \mp && \mp && \mp\end{psmallmatrix}$\\
					&    & or& $\begin{psmallmatrix}
						\alpha, & \beta & {}\lessgtr{} & \gamma & {}\lessgtr{} & \delta, & \alpha, & \beta & {}\lessgtr{} & \gamma & {}\lessgtr{} & \delta \\
						\mp     & \pm   &              & \pm    &              & \pm     & \mp     & \pm   &              & \pm    &              & \pm
					\end{psmallmatrix}$\\
					6 & &   & $\begin{psmallmatrix}
						\alpha, & \beta & {}\gtrless{} & \gamma & {}\lessgtr{} & \delta, & \alpha, & \beta & {}\gtrless{} & \gamma & {}\lessgtr{} & \delta \\
						\pm     & \mp   &              & \mp    &              & \mp     & \pm     & \mp   &              & \mp    &              & \mp
					\end{psmallmatrix}$\\
					7 & &   & $\begin{psmallmatrix}
						\alpha, & \beta, & \gamma, & \delta, & \alpha, & \beta, & \gamma, & \delta \\
						+       & +      & -       & -       & +       & +      & -       & -
					\end{psmallmatrix}$ \\\midrule
					8 & \multirow{5}{*}{$C$} &   & $(++-++-+-)$\\
					9 &  &  & $(--+--+-+)$\\
					10 & &   & $(--++-+-+)$\\
					11 & &   & $(++--+-+-)$\\
					12 & &   & $(--+-++-+)$\\\bottomrule
				\end{tabular}
			}%
		}%
		\hfill
		\raisebox{.5\depth-0.5\height}{%
			\begin{tabular}{@{}c@{}}
				\subcaptionbox{One-dimensional strata\label{tab:1-dim types}}{
					{
						\begin{tabular}{c|cc}
							\toprule
							type  && condition \\\midrule
							13 && $(\alpha, \beta, \gamma, \delta, -\delta, -\gamma, -\beta, -\alpha)$ \\
							14 && $(\alpha, \beta, \gamma, \beta, \alpha, \delta, \epsilon, \zeta)$ \\
							15 && $\begin{psmallmatrix} \alpha,& \beta,& -\alpha,& -\beta,& \alpha,& \beta,& -\alpha,& -\beta  \\ + & + & - & - & + & + & - & - &\end{psmallmatrix}$\\
							\multirow{2}{*}{16} && $\begin{psmallmatrix} \alpha, & \beta &{}\lessgtr{} & \gamma &{}={}&\delta,& \alpha, & \beta &{}\lessgtr{} & \gamma &{}={}&\delta \\ \pm & \mp && \mp && \mp & \pm & \mp && \mp && \mp\end{psmallmatrix}$ \\
							& or & $\begin{psmallmatrix} \alpha, & \beta &{}={} & \gamma &{}\lessgtr{}&\delta,& \alpha, & \beta &{}={} & \gamma &{}\lessgtr{}&\delta \\ \mp & \pm && \pm && \pm & \mp & \pm && \pm && \pm\end{psmallmatrix}$ \\
							17 && $(\alpha, \beta, \gamma, \alpha, \epsilon, -\gamma, -\beta, \epsilon)$ \\
							18 && $(\alpha, \beta, \gamma, \delta, -\alpha, -\beta, -\gamma, -\delta)$ \\
							\bottomrule
						\end{tabular}
					}%
				} \\[24pt]
				\subcaptionbox{Zero-dimensional strata. For examples of these configurations, see \cref{fig:0-dim types}.\label{tab:0-dim types}}{
					\begin{tabular}{ccc}
						\toprule
						intersecting & cardi- & pattern \\
						types & nality & \\\midrule
						3,13,14 & 2 & $(\alpha, -\alpha, \alpha, -\alpha, \alpha, -\alpha, \alpha, -\alpha, )$ \\ 
						13,15 & 4 & $(\alpha, 0, -\alpha, 0, \alpha, 0, -\alpha, 0)$ \\ 
						14,15 & 4 & $(\alpha, \alpha, -\alpha, -\alpha, \alpha, \alpha, -\alpha, -\alpha)$ \\
						13,16 & 8 & $(\alpha, -\alpha, 0, 0, \alpha, -\alpha, 0, 0)$ \\
						14,16 & 8 & $(\alpha, \beta, \beta, \beta, \alpha, \beta, \beta, \beta)$ \\ 
						13,18 & 8 & $(\alpha, -\alpha, \beta, \beta, -\alpha, \alpha, -\beta, -\beta)$\\ 
						14,18 & 8 & $(\alpha,-\beta,0,\beta,-\alpha,\beta,0,-\beta)$ \\ 
						\bottomrule
					\end{tabular}
				}
			\end{tabular}
		}%
	}
	\caption{Dihedral angle patterns of the symmetry types.
		A configuration belongs to a certain symmetry type if its dihedral angle sequence or any circular shift of it satisfies the condition.
		Conditions with signs $\pm$ or $\lessgtr$ are meant as either the condition with all upper signs, or the condition with all lower signs.
		In types 4,5,6 and 16, the magnitude of $\alpha$ compared to $\beta$, $\gamma$ or $\delta$ is not relevant.
	}
	\label{tab:symmetry types}
\end{table}

These symmetry types partition $\Cfg(K)$ into various subspaces, which we analyze separately in the following.
Note that a point can have more than one symmetry types.
It follows from the definition that every point has either type 1 or type 2.

\begin{figure}[tbp]
	\hbadness=10000
	\parfillskip=0pt
	\subcaptionbox*{type  1}{\includegraphics[height=2cm]{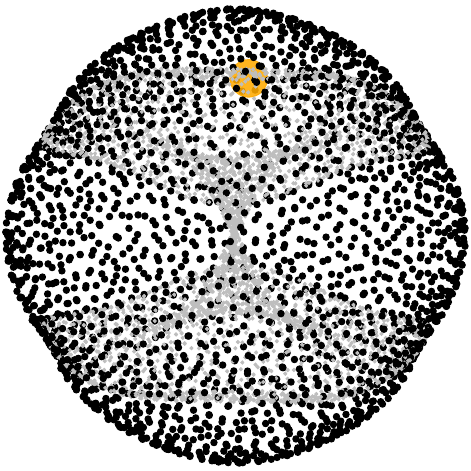}\includegraphics[height=2cm]{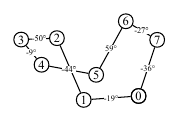}}
	\subcaptionbox*{type  2}{\includegraphics[height=2cm]{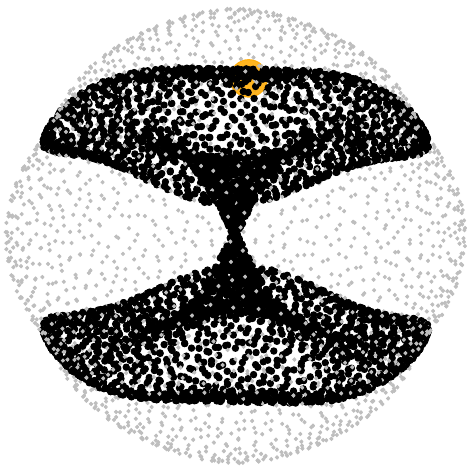}\includegraphics[height=2cm]{new_images/6000/individual_conformations/2335}}
	\subcaptionbox*{type  3}{\includegraphics[height=2cm]{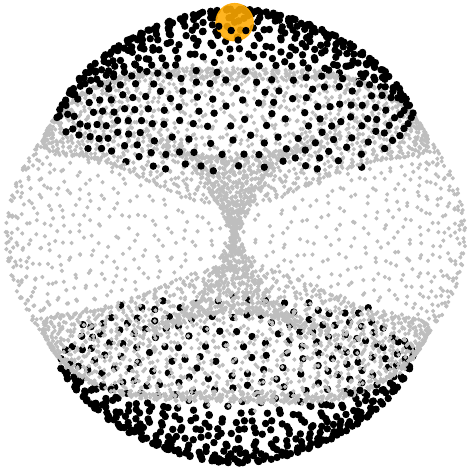}\includegraphics[height=2cm]{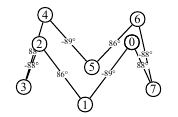}}
	\subcaptionbox*{type  4}{\includegraphics[height=2cm]{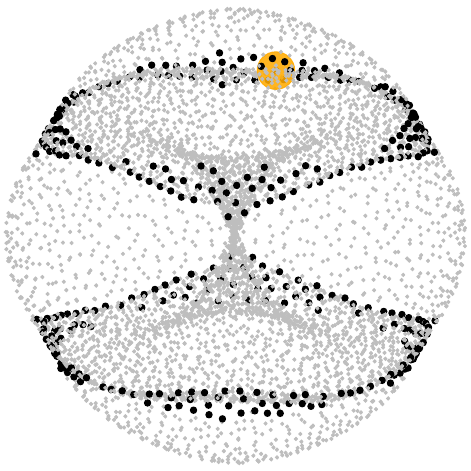}\includegraphics[height=2cm]{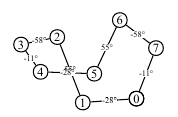}}
	\subcaptionbox*{type  5}{\includegraphics[height=2cm]{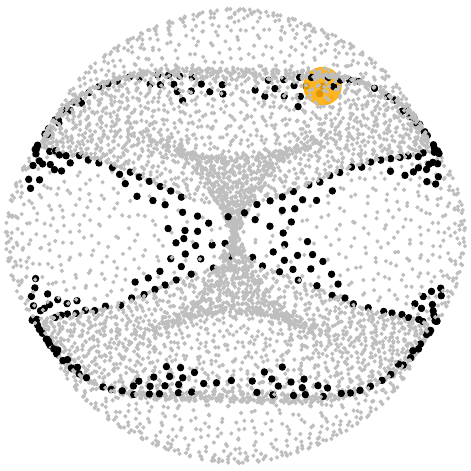}\includegraphics[height=2cm]{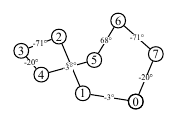}}
	\subcaptionbox*{type  6}{\includegraphics[height=2cm]{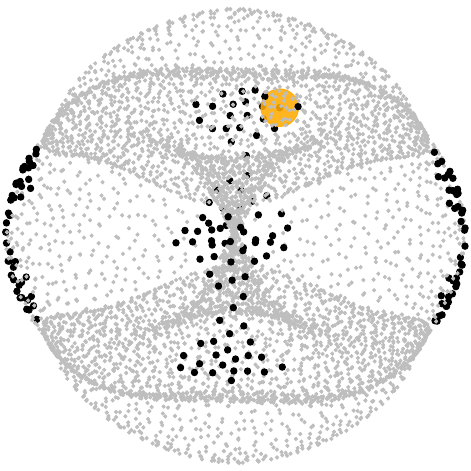}\includegraphics[height=2cm]{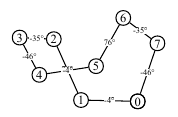}}
	\subcaptionbox*{type  7}{\includegraphics[height=2cm]{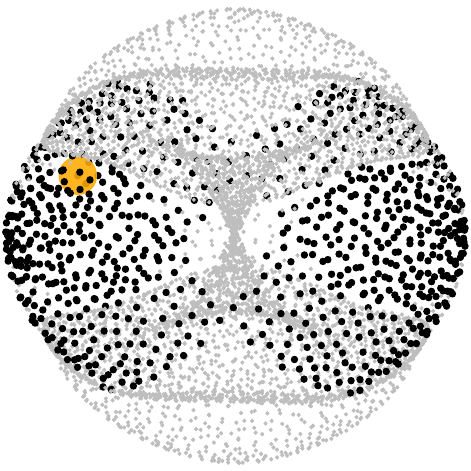}\includegraphics[height=2cm]{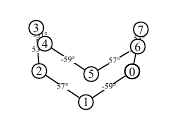}}
	\subcaptionbox*{type  8}{\includegraphics[height=2cm]{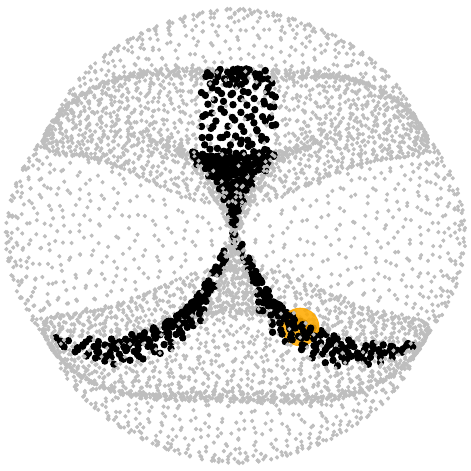}\includegraphics[height=2cm]{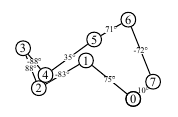}}
	\subcaptionbox*{type  9}{\includegraphics[height=2cm]{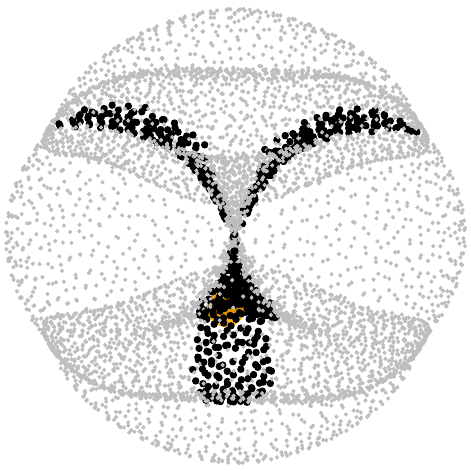}\includegraphics[height=2cm]{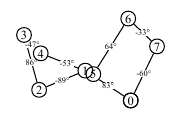}}
	\subcaptionbox*{type 10}{\includegraphics[height=2cm]{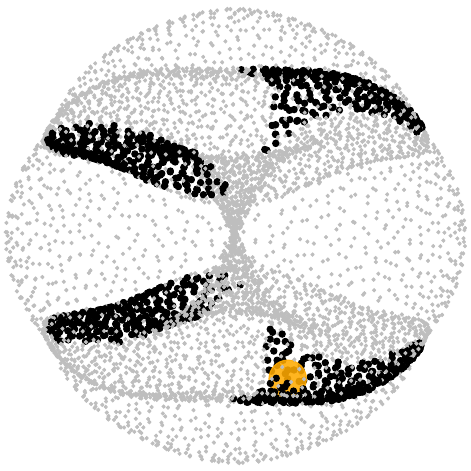}\includegraphics[height=2cm]{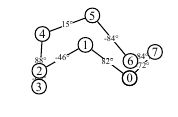}}
	\subcaptionbox*{type 11}{\includegraphics[height=2cm]{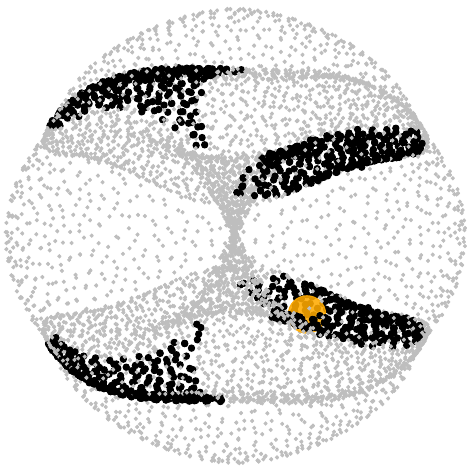}\includegraphics[height=2cm]{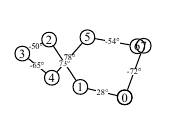}}
	\subcaptionbox*{type 12}{\includegraphics[height=2cm]{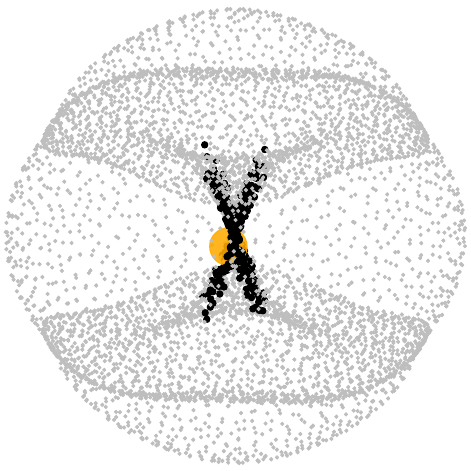}\includegraphics[height=2cm]{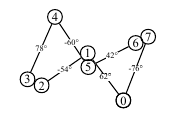}}
	\subcaptionbox*{type 13\label{fig:t13}}{\includegraphics[height=2cm]{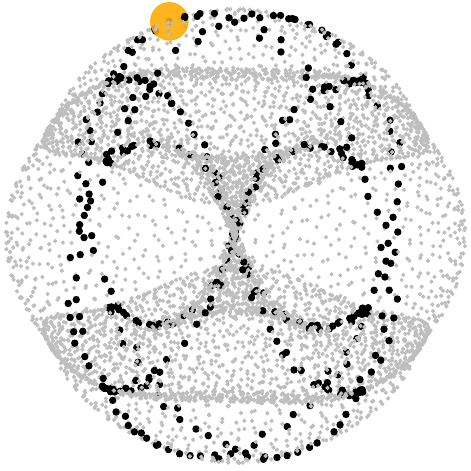}\includegraphics[height=2cm]{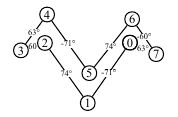}}
	\subcaptionbox*{type 14}{\includegraphics[height=2cm]{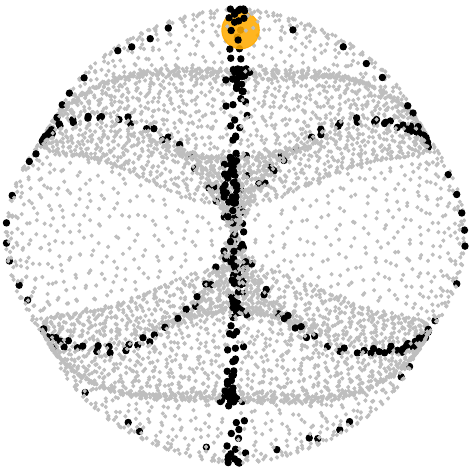}\includegraphics[height=2cm]{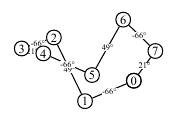}}
	\caption{Location of symmetry types in the Isomap projection w.r.t.\ $d_\Vert$ from \cref{fig:isomap subspaces}.
		The right picture shows a typical standard representatives, whose location in the projection on the left is highlighted.
		See also \cref{fig:symmetry-types angular} for a different projection of the symmetry types on the “hourglass” subspace $A$.}
	\label{fig:symmetry-types}
	\vspace{-13pt}
\end{figure}

\begin{figure}[tbp]
	\samenumber
	\parfillskip=0pt
	\hbadness=10000
	\subcaptionbox*{type 15}{\includegraphics[height=2cm]{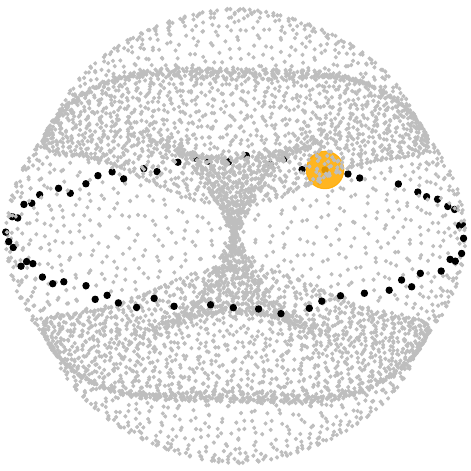}\includegraphics[height=2cm]{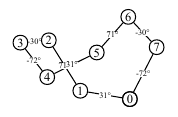}}
	\subcaptionbox*{type 16}{\includegraphics[height=2cm]{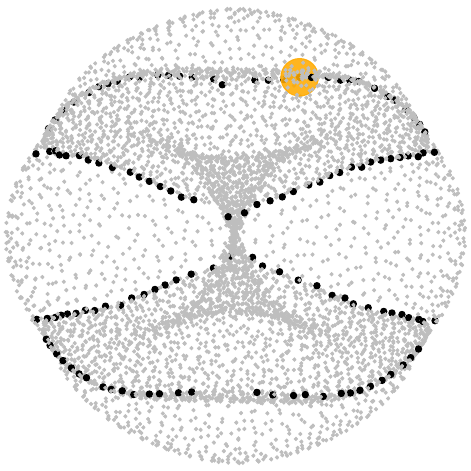}\includegraphics[height=2cm]{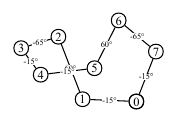}}
	\subcaptionbox*{type 17}{\includegraphics[height=2cm]{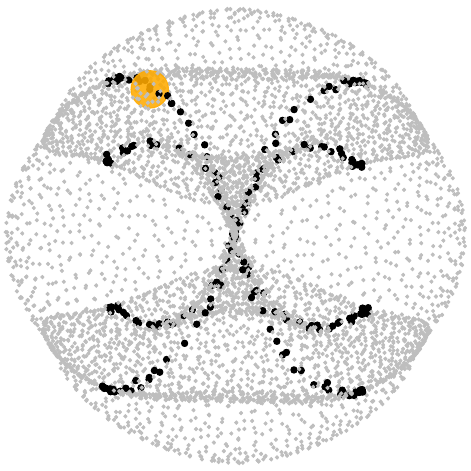}\includegraphics[height=2cm]{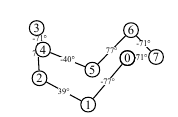}}
	\subcaptionbox*{type 18}{\includegraphics[height=2cm]{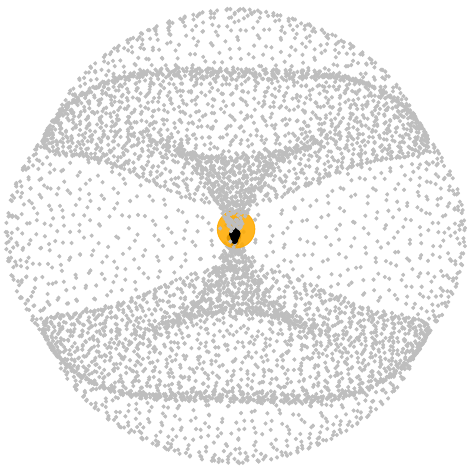}\includegraphics[height=2cm]{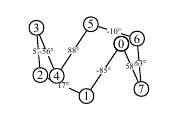}}
	\caption{Location of symmetry types (continued).}
\end{figure}

\begin{figure}[tbp]
	\hbadness=10000
	\parfillskip=0pt
	\subcaptionbox*{type  8}{\includegraphics[height=2cm]{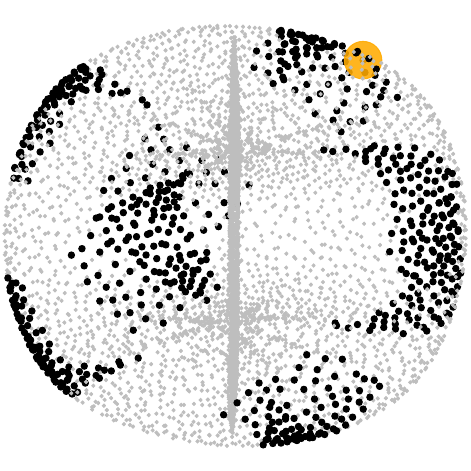}}
	\subcaptionbox*{type  9}{\includegraphics[height=2cm]{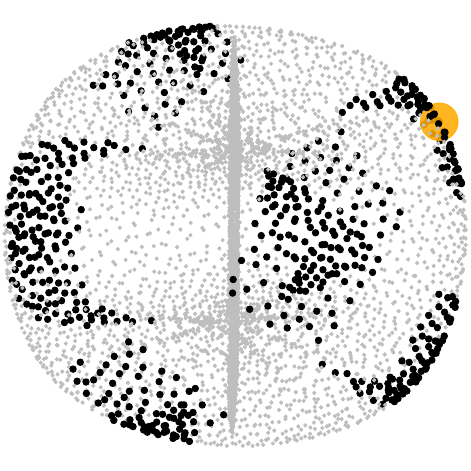}}
	\subcaptionbox*{type 10}{\includegraphics[height=2cm]{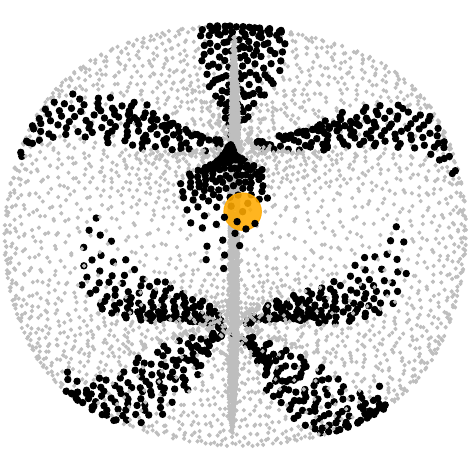}}
	\subcaptionbox*{type 11}{\includegraphics[height=2cm]{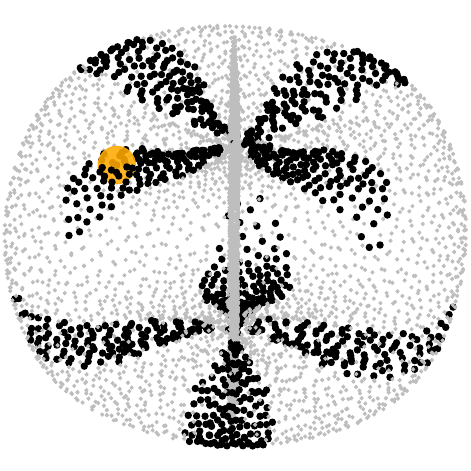}}
	\subcaptionbox*{type 12}{\includegraphics[height=2cm]{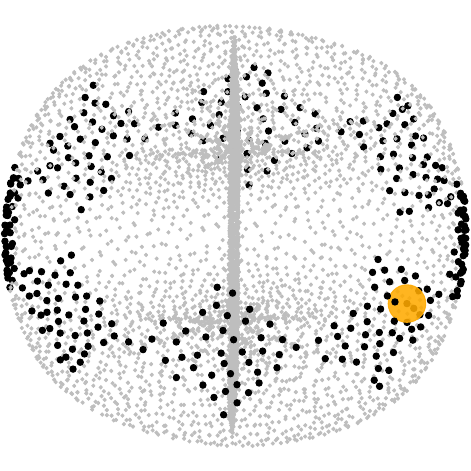}}
	\subcaptionbox*{type 14}{\includegraphics[height=2cm]{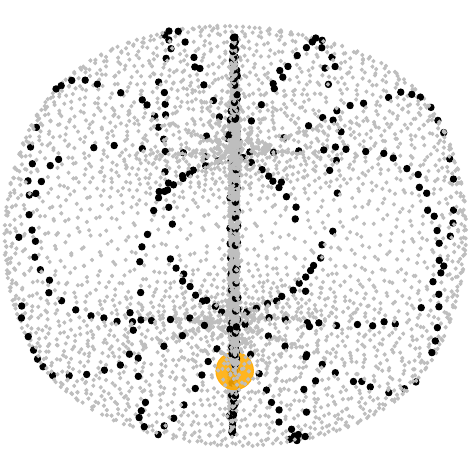}}
	\subcaptionbox*{type 13 and 17\label{fig:t13a}}{\includegraphics[height=2cm]{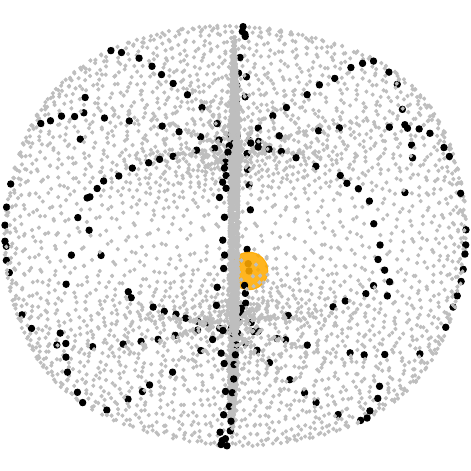}}
	\subcaptionbox*{type 18}{\includegraphics[height=2cm]{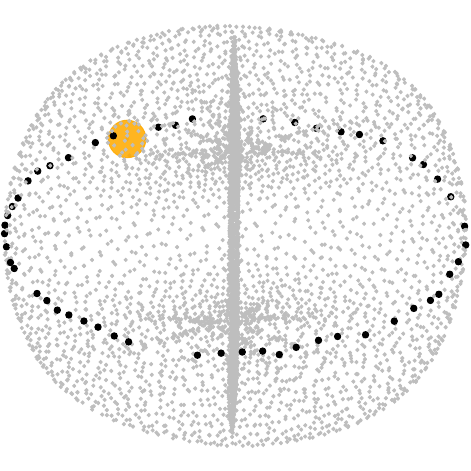}}
	\caption{Location of symmetry types existing on the Klein bottle subspace $C$, in the Isomap projection w.r.t.\ $d_\angle$ from \cref{fig:isomap subspaces}.}
	\label{fig:symmetry-types angular}
	\vspace{0pt}
\end{figure}

\runinhead{Two-dimensional strata}
	First, we say a configuration has symmetry type $1$ if its dihedral angles form a sequence of the form $(\alpha, \beta, \gamma, \delta, \alpha, \beta, \gamma, \delta)$,
	and we say that it has type $2$ otherwise; see \cref{tab:2-dim types}.
	A configuration of type 1 is said to have also symmetry type 3 if its dihedral angles have signs of the form $(+,-,+,-,\dotsc)$ or $(-,+,-,+,\dotsc)$, with $+$ standing for non-negative and $-$ standing for non-positive,
	and of symmetry type 4 if any circular shift $(\alpha, \beta, \gamma, \delta,\dotsc)$ of its dihedral angle sequence
	either has signs $(+,-,-,-,\dotsc)$ and satisfies $\beta \leq \gamma \geq \delta$,
	or has signs $(-,+,+,+,\dotsc)$ and satisfies $\beta \geq \gamma \leq \delta$.

	We let $A \subseteq \Cfg(K)$ be the subspace of all configurations that have symmetry type 3 or 4;
	we let $B \subseteq \Cfg(K)$ be the subspace of all configurations that have symmetry type 1 but not 3 or 4,
	and we let $C \subseteq \Cfg(K)$ be the closure of the subspace of all configurations that have symmetry type 2.
	Then $A$, $B$ and $C$ are two-dimensional subspaces of $\Cfg(K)$.
	Because type 2 is the complement of type 1, they cover $\Cfg(K)$.
	Furthermore, their interiors are pairwise disjoint.
	\Cref{tab:2-dim types} lists further symmetry types (5--12) that partition $\Cfg(K)$ even further.
	All these types define two-dimensional subspaces of $\Cfg(K)$.

	\Cref{fig:symmetry-types,fig:symmetry-types angular} show a configuration of each of these (and the following) symmetry type,
	together with an indication where these symmetry types occur in the Isomap projection
	from \cref{fig:isomap subspaces}.

\runinhead{One-dimensional strata}
	Additionally, we define the symmetry types defined by the conditions in \cref{tab:1-dim types}.
	All of these define one-dimensional subspaces of $\Cfg(K)$.
	Symmetry type 16 is the intersection of types 4 and 5 and will turn out later to be the pairwise intersection of $A$, $B$ and $C$.
	Because of their appearance in \cref{fig:symmetry-types,fig:symmetry-types angular}, we will refer to the types 13, 14 and 17 as \emph{meridians}, and
	to type 15 (which lies in $B$) and type 18 (which lies in $C$) as \emph{equators}.
	Of course, intersections of the closures of the one-dimensional strata defines an additional abundance of one-dimensional strata;
	however, we refrain from giving names to all of them.

\begin{figure}[tbp]
	\makebox[\linewidth][c]{
		\subcaptionbox{3,13,14}{\includegraphics[height=2cm]{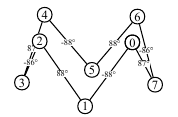}}
		\subcaptionbox{13,15}{\includegraphics[height=2cm]{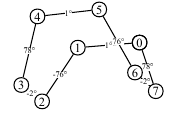}}
		\subcaptionbox{14,15}{\includegraphics[height=2cm]{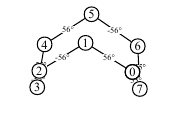}}
	}
	\makebox[\linewidth][c]{
		\subcaptionbox{13,16}{\includegraphics[height=2cm]{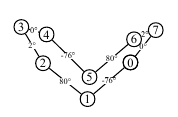}}
		\subcaptionbox{14,16}{\includegraphics[height=2cm]{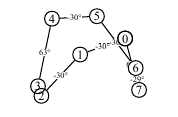}}
		\subcaptionbox{13,18}{\includegraphics[height=2cm]{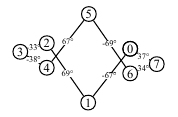}}
		\subcaptionbox{14,18}{\includegraphics[height=2cm]{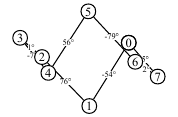}}
	}
	\caption{Examples for the zero-dimensional strata from \cref{tab:0-dim types}.}
	\label{fig:0-dim types}
\end{figure}
\runinhead{Zero-dimensional strata}
	Finally, the intersections of the (closures of) strata listed in \cref{tab:0-dim types} contain only finitely many elements that are particularly symmetric; see \cref{fig:0-dim types}.

\medskip
For $S \subseteq \mathbf{N}$, we denote by $\Cfg(K, S)$ the subspace of $\Cfg(K)$ of configurations
that have a symmetry type in $S$.
In the following, we analyze the subspaces $A = \Cfg(K, \{3,4\})$, $B=\Cfg(K, \{5,6,7\})$ and $C = \Cfg(K, \{2\}) = \Cfg(K) \setminus \Cfg(K, \{1\})$ separately.

Note:
If $\Cfg(K, S)$ defines a zero- or one-dimensional stratum,
we cannot expect to find configurations for which the defining conditions hold on the nose.
Instead, we test the equations defining these symmetry types up to a certain tolerance
($\pm 5^\circ$ for one-dimensional and $\pm 2.5^\circ$ for two-dimensional strata).
In particular, points lying close to the boundary of a component will also be assigned the symmetry type
of the neighboring component.

\subsection{\texorpdfstring{The spherical subspace $A \cup B$}{The spherical subspace}}
\begin{figure}[tbp]
	\centering
	\includegraphics[scale=0.5]{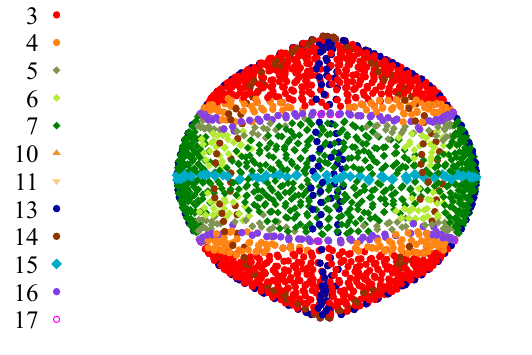}
	\includegraphics{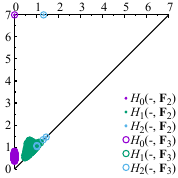}
	\caption{Isomap projection and persistence diagram (w.r.t.\ $d_\angle$) the component $A \cup B$,
		defined as the set of points with symmetry types 3--7.
		Note that a point in $A \cup B$ can \emph{additionally} have one ore more of the symmetry types 13--18
		that define one-dimensional strata.
		Furthermore, because we test the conditions from \cref{tab:symmetry types} only up to a tolerance,
		we also see points that have symmetry types 10 and 11, in addition to one of the types 3--7.
	}
	\label{fig:AB}
\end{figure}

\Cref{fig:AB} shows an Isomap projection of the subspace $A \cup B$ alone, together with the symmetry types of the samples.
These projections suggest that $A \cup B$ is homeomorphic to a 2-sphere $S^2$.
Since there are no visible singularities in the projection, we assume that $A \cup B$ is a surface.
Consider the persistence diagram of $A \cup B$ in \cref{fig:AB}.
It shows that the $\F_2$- and $\F_3$-Betti numbers are (1,0,1) in dimensions (0,1,2), respectively.
In particular, they coincide, hence $A \cup B$ is orientable.
By the classification of surfaces (see \cref{sec:classification of surfaces}), it follows that $A \cup B$ is homeomorphic to $S^2$.

According to \cref{fig:AB}, the sphere
$A \cup B$ decomposes further into the following subspaces:
$A = \Cfg(K,\{3,4\})$ forms two opposing “pole caps” on the sphere,
and $B=\Cfg(K,\{5,6,7\})$ forms a cylinder between the two discs.
The 1-dimensional subspace $\Cfg(K,\{13\})$ forms four “meridians” on the two discs and the cylinder,
which partition $B$ and each of the two connected subspaces of $A$ into four subspaces,
giving the cell structure shown in \cref{fig:cell structure:A,fig:cell structure:B}.
The colors in \cref{fig:cell structure} correspond to the symmetry types in \cref{fig:AB}.

We see that the Isomap projection in \cref{fig:AB} maps $A \cup B$ to a sphere, such that the two connected components of $A$ can be viewed as the “pole caps” of this sphere.
The two poles then are the two configurations $P_1 \coloneqq (\alpha, -\alpha, \alpha, -\alpha,\dotsc)$ and $P_2 \coloneqq (-\alpha, \alpha, -\alpha,\alpha, \dotsc)$, with $\alpha \approx 89^\circ$.
The sphere is divided into two hemispheres $H_1$ and $H_2$, where $H_1$ (resp.\ $H_2$) contains the configurations in $A \cup B$ whose angle sequence $(\alpha_1, \alpha_2,\dotsc)$
satisfies $\alpha_1 + \alpha_3 + \dotsb \geq \alpha_2 + \alpha_4 + \dotsb$ (resp.\ $\leq$).
In particular, $P_i \in H_i$ for $i = 1,2$, and $H_1 \cap H_2 = \Cfg(K, \{15\})$.

\subsection{The Klein bottle subspace \texorpdfstring{$C$}{C}}
This subspace has been referred to in the literature \cite{MartinThompsonEtAl:2010,AdamsMoy:2021} as the “interior” or “hourglass” subspace,
in reference to its appearance as the hourglass inscribed into the sphere
in \cref{fig:isomap subspaces}.
It has been suggested that this subspace can be described as a Klein bottle \cite{MartinThompsonEtAl:2010}.
We provide additional evidence for this observation, extending the results from \cite{AdamsMoy:2021}.

\begin{figure}[tbp]
	\centering
	\includegraphics[scale=0.5]{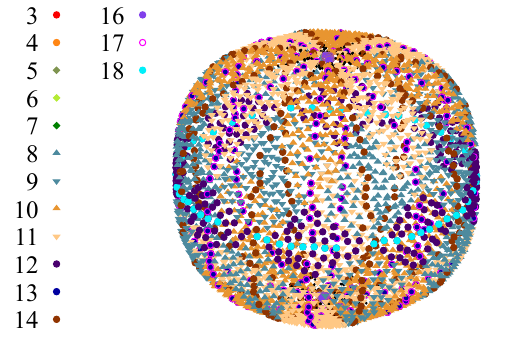}%
	\includegraphics{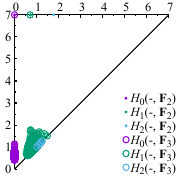}%
	\caption{Isomap projection and persistence diagram (w.r.t.\ $d_\angle$) of $C$.
		See \cref{fig:AB} for an explanation of the appearance of symmetry types 4--7 in this plot.
	}
	\label{fig:klein-bottle:isomap}
	\label{fig:klein-bottle:persistence}
\end{figure}
The Isomap projection in \cref{fig:isomap subspaces} exhibits high metric distortion around the two circles
in which $A \cup B$ and $C$ intersect.
\Cref{fig:klein-bottle:isomap} shows the Isomap projection of $C$ alone.
In this sphere-like projection, this intersection is mapped to the two “poles” of the sphere and is subject to large metric distortion;
i.e., Isomap places points close to each other in the projection that are distant in the actual space.
We suspect that this region is in fact a 1-cell to which the 2-cell forming $C$ is attached by a degree-2-map,
yielding a space that is not isometrically embeddable in $\R^3$.
This would align well with the fact that as an algebraic surface,
$\Cfg(K)$ cannot have one dimensional singularities with an uneven number of leaves attached to it.

Therefore, we propose that $C$ can be represented by the cell complex shown in \cref{fig:cell structure:C}, which refines the complex in \cref{fig:model:C}.
Again, all edges with the same endpoints denote the same 1-cell.
This cell complex is a model for a Klein bottle;
see also \cref{fig:model:C,fig:model:C:2}.
As we will see in \cref{sec:gluing-by-meridians}, the vertices $v_j$ and $w_i$ in \cref{fig:cell structure:C} correspond to the vertices in \cref{fig:cell structure:A,fig:cell structure:B}.

To verify that \cref{fig:cell structure:C} models $C$,
consider the persistence diagram of $C$ in \cref{fig:klein-bottle:persistence}.
It suggests that $C$ has Betti numbers (1, 2, 1) over $\F_2$ and (1,1,0) over $\F_3$ in dimensions 0, 1 and 2, respectively.
Assuming that $C$ is a manifold, the classification of closed surfaces gives that $C$ is a Klein bottle.

\phantomsection
\label{sec:circular coordinates}
To provide additional evidence for this claim,
we have selected only the points from $C$ whose angle sequence $(\alpha_1, \alpha_2,\dotsc)$ satisfies $\alpha_1 + \alpha_3 + \dotsb \geq \alpha_2 + \alpha_4 + \dotsb$.
As for the sphere component $A \cup B$, this condition (and the analogous condition with $\leq$) define two subspaces $M_1$ and $M_2$ (the upper and lower “hemisphere” in \cref{fig:klein-bottle:isomap})
that cover $C$ and satisfy $\Cfg(K, \{18\})$.
This corresponds to cutting the cell complex in \cref{fig:model:C,fig:cell structure:C} along the horizontal line.
If $C$ is a Klein bottle, then $M_1$ and $M_2$ will be two Möbius strips with configurations of type 16 as core curve and type 18 as boundary.

\begin{figure}[tbp]
	\centering
	\subcaptionbox{$M_1$\label{fig:mobius-band}}{%
		\includegraphics[height=3cm]{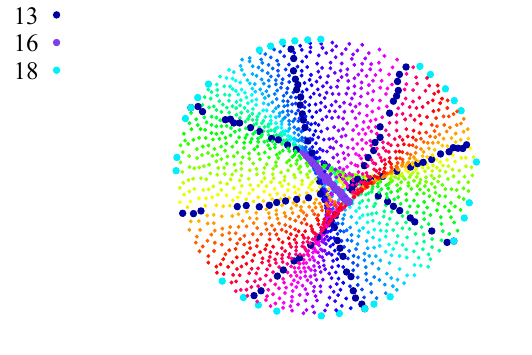}%
		\includegraphics{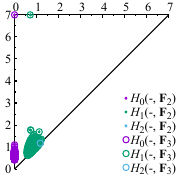}%
	}\hfill
	\subcaptionbox{$M_1 \cap \Cfg(K, \{16\})$\label{fig:mobius-band-boundary}}{%
		\includegraphics[height=3cm, trim={25mm 0 0 0},clip]{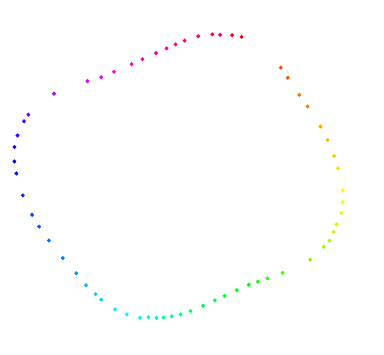}%
	}\hfill
	\subcaptionbox{\label{fig:rainbow}}{\includegraphics[trim={0 0 .5mm 0},clip]{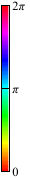}}
	\caption{Isomap projection and persistent homology (w.r.t.\ $d_\angle$) of the Möbius strip $M_1$.
		The rainbow colors correspond to a map $f\colon M_1 \to S^1$ that represents the prominent $1$-cohomology class in the persistence diagram;
		see (\subref{fig:rainbow}).
		Restricting $f$ to the boundary of $M_1$ (the component $\Cfg(K, \{18\}$, which is a circle),
		we see that $f|_{\Cfg(K, \{18\})}$ is a map $S^1 \to S^1$ of mapping degree 2,
		while when restricted to the core curve of $M_1$ (the component $M_1 \cap \Cfg(K, \{16\})$, which also is a circle, see (\subref{fig:mobius-band}))
		yields a map $f|_{M_1 \cap \Cfg(K, \{16\})}\colon S^1 \to S^1$ of degree 1.
		This supports the hypothesis that $M_1$ is a Möbius strip.
	}
\end{figure}

Indeed, the Isomap projection of $M_1$ in \cref{fig:mobius-band} resembles a Plücker conoid, which is a specific immersion of the Möbius strip into $\R^3$.
Furthermore, according to the persistence diagram of $M_1$ in \cref{fig:mobius-band}, $M_1$ has the homology type of a 1-sphere.
Through the bijection $H^1(M_1, \Z) \cong [M_1, S^1]$ between the first integral cohomology of $M_1$ and homotopy classes of maps from $M_1$ to $S^1$,
one can associate a representative map $M_1 \to S^1$ to each class in the persistent cohomology $H^1(\VR_*(M), \Z)$.
For details on the method, see \cite{deSilvaMorozovEtAl:2011}.

The color-coding in \cref{fig:mobius-band} shows such a map $f\colon M_1 \to S^1$ corresponding to the prominent 1-cohomology class in the persistence diagram in \cref{fig:mobius-band}.
We see that $f$ winds two times around the boundary if $M$.
An Isomap projection of the core curve $M_1 \cap \Cfg(K, \{16\})$ of $M_1$, together with the restriction of $f$, is shown in \cref{fig:mobius-band-boundary}.
We see that $f$ winds around this core curve once, as expected.
We therefore argue that $M_1$ (and analogously, $M_2$) are two Möbius strips that are glued together along their boundaries $\Cfg(K, \{18\})$.
Therefore, $C$ is a Klein bottle.

\subsection{Gluing both subspaces}
\label{sec:gluing-by-meridians}
It remains to understand how the two subspaces $A$, $B$ and $C$ and are glued together.
We have already seen that $A$ and $B$ intersect in two circles given by $\Cfg(K, \{16\})$.
In the previous section, we have also seen that the core curves of the two Möbius strips $M_1$, $M_2$ are also $\Cfg(K, \{16\})$.
Thus, $A$, $B$ and $C$ commonly intersect in two circles given by $\Cfg(K, \{16\})$.

To make clear how this gluing works, consider the aforementioned type 13 “meridians” on $\Cfg(K)$; see \cref{fig:t13,fig:t13a}
These meridians bifurcate in the points where $A \cup B$ and $C$ intersect.
We can thus track where these bifurcation points lie in the projections of \cref{fig:isomap subspaces}.

\begin{figure}[tbp]
	\subcaptionbox{$C$ (w.r.t.\ $d_\angle$)              \label{fig:meridians:C}  }[.3\linewidth]{\centering\includegraphics[height=3cm]{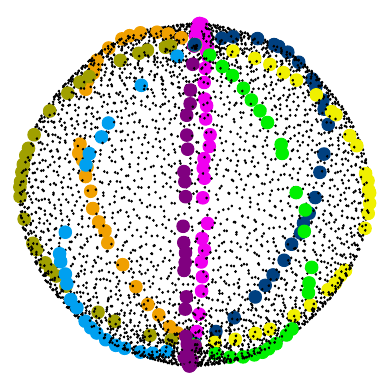}}\hfill
	\subcaptionbox{$\Cfg(K)$ (w.r.t.\ $d_\Vert$)         \label{fig:meridians:ABC}}[.3\linewidth]{\centering\includegraphics[height=3cm]{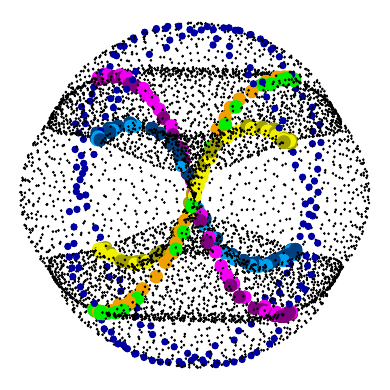}}\hfill
	\subcaptionbox{$\Cfg(K, \{13\})$ (w.r.t.\ $d_\Vert$) \label{fig:meridians:13} }[.3\linewidth]{\centering\includegraphics[height=3cm]{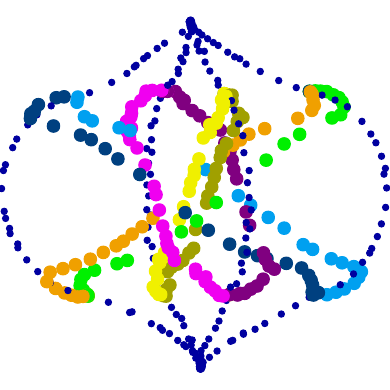}}
	\caption{
		$\Cfg(K, \{13\})$ endows $C$ with a system of “meridians”; cf.\ \cref{fig:t13,fig:t13a}.
		The pictures show the Isomap projections from \cref{fig:isomap subspaces}, where the individual meridians of $C$ have been equipped with colors,
		which allows to understand how the two components $C$ (\subref{fig:meridians:C} and the “hourglass” in \subref{fig:meridians:ABC})
		and $A\cup B$ (the sphere in \cref{fig:meridians:ABC}) are glued together.
	}
	\label{fig:meridians}
\end{figure}

To this end, we use the following observation:
consider the points in $\Cfg(K, \{13\}) \cap C$.
A point of type 13 has the dihedral angle pattern
$(\alpha, \beta, \gamma, \delta, -\delta, -\gamma, -\beta, -\alpha)$
or one of its circular shifts $(\alpha, \beta,\gamma, \delta,\dotsc)$, $(\beta, \gamma,\delta,-\delta,\dotsc),\dotsc$.
There are four distinct circular shifts before one gets back the original pattern.
If one imposes the additional condition $\delta \geq 0$, then one gets eight distinguishable shifts of this pattern,
In \cref{fig:meridians}, we have highlighted these eight distinct shifts of this pattern on $\Cfg(K, \{13\}) \cap C$ in a different color.

In the Isomap projection of $C$ alone (\cref{fig:meridians:C}), we see that this approach
indeed distinguishes all the eight meridians on the Klein bottle component.
When applying the same color-scheme in the Isomap projection of the entire space (see \cref{fig:meridians:ABC})
in which $C$ is projected to the “hourglass” sitting inside the sphere (which represents $A \cup B$),
we see that each two antipodal meridians of $C$ are projected to the same image in \cref{fig:meridians:ABC}.
This happens because samples of $A \cup B$ have a dihedral angle pattern of the form $(\alpha, \beta, \gamma, \delta, \alpha, \beta, \gamma, \delta)$
and therefore only have four distinct circular shifts.
Furthermore, we see that any pair of antipodal meridians on $C$ (e.g., the two blue ones)
are attached to two antipodal meridians on $A \cup B$; see \cref{fig:meridians:ABC};
see also the Isomap projection of $\Cfg(K, \{13\})$ alone in \cref{fig:meridians:13}.

\begin{figure}[tbp]
	\centering
	\subcaptionbox{subspace $A$\label{fig:cell structure:A}}{
		\begin{tikzpicture}
			\node (v0) at (0,0) {$v_0$};
			\node (w0) at (3,0) {$w_0$};
			\foreach[count=\i] \phi in {180, 270, 0, 90}{
				\node                (v\i) at (\phi:1) {$v_\i$};
				\node[shift={(3,0)}] (w\i) at (\phi:1) {$w_\i$};
				\draw (v0) edge[halo, t13, "{\ifthenelse{\i=4}{13}{}}"' {color=13, ll}] (v\i);
				\draw (w0) edge[halo, t13, "{\ifthenelse{\i=4}{13}{}}"' {color=13, ll}] (w\i);
			}
			\foreach[remember=\i as \j (initially 1), remember=\out as \in (initially 180)] \i/\out in {2/270, 3/0, 4/90, 1/180} {
				\path (v\j) to[out=\out, in=\in] coordinate (g\i) (v\i);
				\path (w\j) to[out=\out, in=\in] coordinate (h\i) (w\i);
				\draw (v0) edge[t14, halo, shorten >=2pt, "{\ifthenelse{\j=4}{14}{}}"' {color=14, ll}] (g\i);
				\draw (w0) edge[t14, halo, shorten >=2pt, "{\ifthenelse{\j=4}{14}{}}"' {color=14, ll}] (h\i);
				\path (v\j) edge[t16, halo, <-, out=\out, in=\in, "{\ifthenelse{\j=4}{16}{}}"' {color=16, ll}] coordinate (g\i) (v\i);
				\path (w\j) edge[t16, halo, <-, out=\out, in=\in, "{\ifthenelse{\j=4}{16}{}}"' {color=16, ll}] coordinate (h\i) (w\i);
			}
			\begin{scope}[on background layer]
				\foreach \a in {v, w} {
					\fill[t4]
						(\a 1.south) to[out=-90,in=180] (\a 2.west)  to (\a 2.center) to (\a1.center) to (\a1.south)
						(\a 2.east)  to[out=  0,in=-90] (\a 3.south) to (\a 3.center) to (\a2.center) to (\a2.east)
						(\a 3.north) to[out= 90,in=  0] (\a 4.east)  to (\a 4.center) to (\a3.center) to (\a3.north)
						(\a 4.west)  to[out=180,in= 90] (\a 1.north) to (\a 1.center)  to (\a 4.center) to (\a4.west);
					\fill[t3] (\a1.center) -- (\a2.center) -- (\a3.center) -- (\a4.center) -- (\a1.center);
				}
				\node[color=4, ll, fill=none] at (-0.75, -0.45) {4};
				\node[color=3, ll, fill=none] at (-0.5, -0.2) {3};
			\end{scope}
		\end{tikzpicture}
	}
	\hfill
	\subcaptionbox{subspace $B$\label{fig:cell structure:B}}{
		\begin{tikzpicture}
			\foreach[count=\i] \j in {1,2,3,4,1}{
				\node (v\i) at (\i, 1) {$v_\j$};
				\node (w\i) at (\i,-1) {$w_\j$};
			}
			\foreach \i in {2,3,4}{
				\draw[halo, t13] (v\i) -- (w\i);
			}
			\draw[halo, t13, ->] (v1) to["$c$"' color=13, "13"' {color=13, ll, near start}] coordinate (c1) (w1);
			\draw[halo, t13, ->] (v5) to["$c$"  color=13] coordinate (c2) (w5);
			\draw[t15, halo] ([xshift=2pt]c1) to["15" {color=15, ll, pos=0.015}] ([xshift=-2pt]c2);
			\foreach [remember=\i as \j (initially 1)] \i in {2,...,5}{
				\draw[<-, halo, t16] (v\i) to["{\ifthenelse{\j=1}{16}{}}"' {color=16, ll}] coordinate (g\j) (v\j);
				\draw[<-, halo, t16] (w\i) to["{\ifthenelse{\j=1}{16}{}}"  {color=16, ll}] coordinate (h\j) (w\j);
				\draw ([yshift=-2pt]g\j) edge[t14, halo, "{\ifthenelse{\j=1}{14}{}}"' {color=14, ll, near start} ] ([yshift=2pt]h\j);
				\coordinate (u\j) at ($(0.5,0)+(\j,0)$);
				\path (u\j) to coordinate (s\j-1) (v\j)
				      (u\j) to coordinate (s\j-2) (v\i)
				      (u\j) to coordinate (s\j-3) (w\j)
				      (u\j) to coordinate (s\j-4) (w\i);
				\scoped[on background layer] {
					\fill[t5] (v\j) -- (s\j-1) -- (g\j) -- (v\j)
					          (v\i) -- (s\j-2) -- (g\j) -- (v\i)
					          (w\j) -- (s\j-3) -- (h\j) -- (w\j)
					          (w\i) -- (s\j-4) -- (h\j) -- (w\i);
					\fill[t6] (u\j) -- (s\j-1) -- (g\j) -- (s\j-2)
					          (u\j) -- (s\j-3) -- (h\j) -- (s\j-4);
					\fill[t7] (v\j.center) -- (w\i.center) -- (v\i.center) -- (w\j.center) -- (v\j.center);
				}
			}
			\node[color=5!75!black, ll, fill=none] at (1.25,-0.75) {5};
			\node[color=6!75!black, ll, fill=none] at (1.4,-0.5) {6};
			\node[color=7, ll, fill=none] at (1.2,-0.25) {7};
		\end{tikzpicture}
	}
	\subcaptionbox{subspace $C$\label{fig:cell structure:C}}{
		\begin{tikzpicture}
			\foreach[count=\j] \i/\k/\l in {1/1/3, 2/2/4, 3/3/1, 4/4/2, 1/5/3, 2/6/4, 3/7/1, 4/8/2, 1/1/3} {
				\node (v\j) at (\j, 1) {$v_\i$};
				\node (w\j) at (\j,-1) {$w_\l$};
			}
			\foreach \i in {1,...,8}{
				\draw[halo, t17] (v\i) to coordinate[pos=.25] (r\i-1) coordinate[pos=.75] (r\i-2) (w\i);
			}
			\path[t17, ->] (v1) to["$d$"' color=13, "13"' {color=13, ll, near start}] coordinate (d1) (w1);
			\path[t17, ->] (v9) to["$d$"  color=13] coordinate (d2) coordinate[pos=.25] (r9-1) coordinate[pos=.75] (r9-2) (w9);
			\draw[t18, halo] ([xshift=2pt]d1) to["18" {color=18, ll, pos=0.015}] ([xshift=-2pt]d2);
			\foreach [evaluate=\i as \j using int(\i+1)] \i in {1,...,8}{
				\draw[->, t16, halo] (v\i) to["{\ifthenelse{\i=1}{16}{}}"  {color=16, ll}] coordinate (g\j) (v\j);
				\draw[->, t16, halo] (w\i) to["{\ifthenelse{\i=1}{16}{}}"' {color=16, ll}] coordinate (h\j) (w\j);
				\draw ([yshift=-2pt]g\j) edge[t14, halo, "{\ifthenelse{\i=1}{14}{}}"' {color=14, ll, near start}] coordinate[pos=.125](g\j') coordinate[pos=.875] (h\j') coordinate (s\i) ([yshift=2pt]h\j);
				\scoped[on background layer]{
					\fill[t12] (r\i-1) -- (r\j-2) -- (r\j-1) -- (r\i-2);
					\fill[t8] (r\i-1) -- (r\j-2) -- (h\j') -- (r\i-2) -- (r\j-1) -- (g\j') -- cycle;
					\fill[t10] (r\i-1) -- (g\j') -- (r\j-1) -- (v\j.center) -- (g\j') -- (v\i.center) -- cycle
					(r\i-2) -- (h\j') -- (r\j-2) -- (w\j.center) -- (h\j') -- (w\i.center) -- cycle;
				}
			}
			\node[color=10!75!black, ll, fill=none] at (1.75,-0.8) {10,11};
			\node[color=8!50!black, ll, fill=none] at (1.3,-0.5) {8,9};
			\node[color=12, ll, fill=none] at (1.1,-0.2) {12};
		\end{tikzpicture}
	}
	\caption{%
		Proposed cell complex representing $\Cfg(K)$.
		The colors represent the symmetry types (small numbers, same colors as in \cref{fig:AB,fig:klein-bottle:isomap}).
		The domains of the different symmetry types is obtained from the Isomap projections.
		For clarity, not all symmetry types listed in \cref{tab:symmetry types} are explicitly drawn or labeled.
		Areas left white have only symmetry type 2 assigned.
		This cell complex is a refinement of \cref{fig:cell complex}.
	}
	\label{fig:cell structure}
\end{figure}
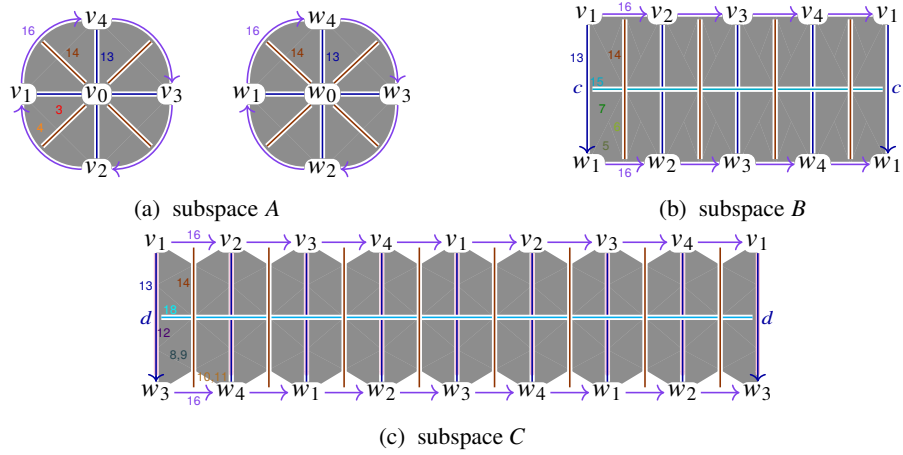

\medskip
From this, we obtain the gluing shown in \cref{fig:cell structure}.
As mentioned before, this cell complex has the $\mathbf{F}_2$- and $\mathbf{F}_3$-Betti numbers (1,1,2) in dimensions (0,1,2), respectively,
which matches the observed persistent homology from \cref{fig:persistence:config-space}.

\section{A cell structure for the quotient space \texorpdfstring{$\Cfg(K)/C_8$}{Cfg(K)/C8}}
Recall that $\Cfg(K)$ is the \emph{labeled} configuration space.
We are also interested in a cell structure of the \emph{unlabeled} configuration space.
In this space, two configurations are considered equivalent if they differ only by a cyclic renumbering of the vertices,
or a reversal of the numbering.
Mathematically, this corresponds to taking the quotient space $\Cfg(K)/D_8$
with respect to the \emph{dihedral group} $D_8$ (the symmetry group of the regular octagon),
which is the semidirect product $D_8 = C_8 \rtimes C_2$ of two cyclic groups,
where $C_8$ acts by cyclic shifting of the numbering, and $C_2$ acts by reversing the numbering.
Thus, our goal is to find a cell structure of $\Cfg(K)/D_8$, and to describe the quotient map $\Cfg(K) \to \Cfg(K)/D_8$.

The quotient $\Cfg(K) \to \Cfg(K)/D_8$ factors through the quotient $\Cfg(K)/C_8$ with respect to the normal subgroup $C_8 \triangleleft D_8$.
Therefore, in this section, we investigate the quotient space $\Cfg(K)/C_8$ as an intermediate step,
together with the quotient map $\Cfg(K) \to \Cfg(K)/C_8$.
We then descibe $\Cfg(K)/D_8$ and the quotient map $\Cfg(K)/C_8 \to \Cfg(K)/D_8$ in \cref{sec:Configuration-space-mod-D8}.

The following definition turns $\Cfg(K)/C_8$ into a metric space.

\begin{definition*}
	Let $X$ be a metric space and $Y=X/{\sim}$
	be a quotient of $X$ w.r.t.\ an equivalence relation $\sim$.
	The \emph{quotient (pseudo-)metric} is the finest
	metric $d_Y$ on $Y$	such that $d_Y([x], [y])\leq d_X(x,y)$ for all $x,y\in X$.
	If all equivalence classes of $X$ w.r.t.\ $\sim$ are finite, $d_Y$ is a metric.
\end{definition*}

\subsection{Actions}
\label{sec:C8 action-on-configurations}
Let $t$ denote the generator of $C_8$, which acts on the vertices by shifting the indexing by one.
On the torsion angles $(\sigma_i(\bar{x}))_i$ of a configuration $\bar{x}\in \Cfg(K)$, the generator $t$ acts by
\[
	t(\sigma_0(\bar{x}),\dotsc,\sigma_7(\bar{x})) = (\sigma_1(\bar{x}),\dotsc,\sigma_7(\bar{x}),\sigma_0(\bar{x})).
\]
The angular metric $d_\angle$ on $\Cfg(K)/C_8$ then is given by
\[
	d_\angle([\bar{x}], [\bar{y}]) = \min \{ d_\angle\bigl(\sigma(\bar{x}), t^n \sigma(\bar{y})\bigr) ; n=0,\dotsc,7 \}.
\]
For the Euclidean metric, the situation is a little more intricate.
Recall from \cref{sec:defining-metrics} that we have defined the Euclidean metric $d_\Vert$ on $\Cfg(K)$ in terms of standard realizations $x^\mathrm{std}$
that satisfy the Eckart condition.
On standard realizations, $t$ acts by
\begin{equation*}
	(x_0,\dotsc,x_7) \mapsto (T x_1,\dotsc, Tx_7, Tx_0)
	\qquad
	\text{for}
	\quad
	T=\begin{psmallmatrix}
		\sqrt 2/2 & \sqrt 2/2 & 0\\
		-\sqrt 2/2& \sqrt 2/2 & 0\\
		0 & 0 & 1
	\end{psmallmatrix},
\end{equation*}
i.e., by cyclic reindexing and a rotation by $\frac{\pi}{4}$ about the $z$-axis.
One may check that this action indeed sends standard realizations to standard realizations.
Thus, the distance of two equivalence classes $[x^\text{std}]$ and $[y^\text{std}]$ in $\Cfg(K)/C_8$ is
\[
	d_\Vert([x^\text{std}], [y^\text{std}])=\min\{ d_\Vert(t^n x, y) ; n = 1,\dotsc,7\}.
\]

\subsection{Equivalence classes of specific configurations}
We start our analysis with the zero dimensional strata from \cref{tab:0-dim types}.
There are two configurations jointly of types 1, 6 and 8 (pattern: $(\alpha,-\alpha,\alpha,-\alpha,\alpha,-\alpha,\alpha,-\alpha)$), corresponding to the two possible signs for $\alpha$.
Each of the two configuration sits in the middle of the two connected components of subspace $A$ (see \cref{fig:isomap subspaces}).
These belong to the same equivalence class mod $C_8$, which we denote by $v_0'$.

Analogously, all four configurations jointly of type 13 and 15 (pattern: $(\alpha,\allowbreak 0,\allowbreak -\alpha,\allowbreak 0,\allowbreak \alpha,\allowbreak 0,\allowbreak -\alpha,\allowbreak 0)$)
(each one sits in the middle of one face of $B$ in \cref{fig:cell structure})
belong to the same $C_8$-orbit, and similarly for the eight configurations of type 13 and 18 (pattern: $(\alpha, -\alpha, \beta, \beta, -\alpha, \alpha, -\beta, -\beta)$).
This suggests that each of the $2$-cells in \cref{fig:cell structure:B,fig:cell structure:C} is a fundamental domain for the $C_8$-action on $B$ and $C$, respectively.
As a justification, recall that the four 2-cells of $A$ and $B$ (see \cref{fig:cell structure})
correspond to the four possible shifts of the dihedral angle pattern $(\alpha,\beta,\gamma,\delta,\alpha,\beta,\gamma,\delta)$.
An analogous justification holds for configurations in $C$.

The four configurations jointly of type 13 and 16 (pattern: $(\alpha,\allowbreak -\alpha,\allowbreak 0,\allowbreak 0,\allowbreak \alpha,\allowbreak -\alpha,\allowbreak 0,\allowbreak 0,\allowbreak )$) that lie in the intersections of $A$, $B$ and $C$
belong to \emph{two} $C_8$ orbits, one for each sign of $\alpha$.
We denote these by $v'_1$ and $v'_2$.

\subsection{Cell structure}
Based on these observations, we arrive at the following description of the cell structure of $\Cfg(K)/C_8$ and the quotient map $\Cfg(K) \to \Cfg(K)/C_8$:
In \cref{fig:mod-C8:projection}, we have replaced the labels of the vertices
by the $C_8$-equivalence classes they belong to.
Edges with identical labels are identified.
The resulting quotient cell complex is shown in \cref{fig:cell structure C8}.
The subcomplex $A' \coloneqq A/C_8$ is a disc; the subcomplex $B' \coloneqq B/C_8$ and $C' \coloneqq C/C_8$ both are isomorphic to a Möbius strip $\mathit{MS}$.
This cell complex represents a gluing $D^2\cup_{S_1}\mathit{MS}\cup_{S_1}\mathit{MS}$
of a disc with two Möbius strips along their boundaries.
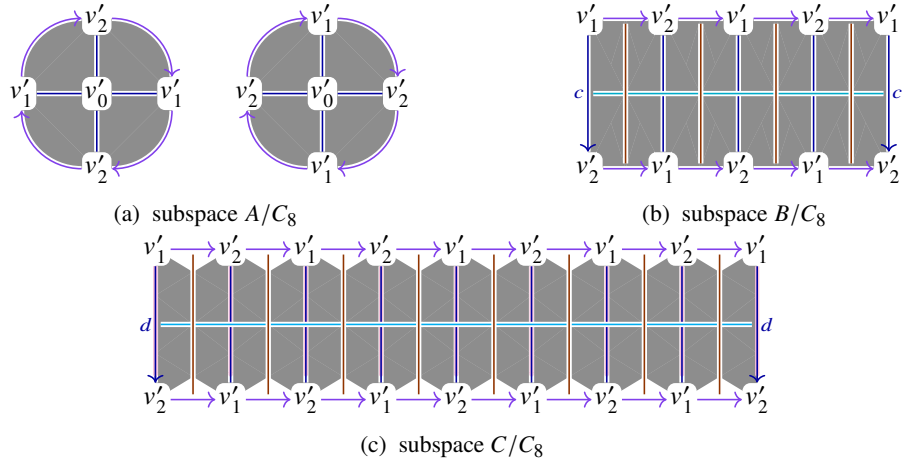
\begin{figure}[tbp]
	\centering
	\subcaptionbox{subspace $A/C_8$}{
		\begin{tikzpicture}
			\node (v0) at (0,0) {$v'_0$};
			\node (w0) at (3,0) {$v'_0$};
			\foreach[count=\i, evaluate={\i as \j using int(mod(\i+1,2)+1)}, evaluate={\i as \k using int(mod(\i+2,2)+1)}] \phi in {180, 270, 0, 90}{
				\node                (v\i) at (\phi:1) {$v'_\j$};
				\node[shift={(3,0)}] (w\i) at (\phi:1) {$v'_\k$};
				\draw[halo, t13] (v0) -- (v\i) (w0) -- (w\i);
			}
			\foreach[remember=\i as \j (initially 1), remember=\out as \in (initially 180)] \i/\out in {2/270, 3/0, 4/90, 1/180} {
				\draw[<-, halo, t16] (v\j) to[out=\out, in=\in] (v\i);
				\draw[<-, halo, t16] (w\j) to[out=\out, in=\in] (w\i);
			}
			\begin{scope}[on background layer]
				\foreach \a in {v, w} {
					\fill[t4]
					(\a 1.south) to[out=-90,in=180] (\a 2.west)  to (\a 2.center) to (\a1.center) to (\a1.south)
					(\a 2.east)  to[out=  0,in=-90] (\a 3.south) to (\a 3.center) to (\a2.center) to (\a2.east)
					(\a 3.north) to[out= 90,in=  0] (\a 4.east)  to (\a 4.center) to (\a3.center) to (\a3.north)
					(\a 4.west)  to[out=180,in= 90] (\a 1.north) to (\a 1.center)  to (\a 4.center) to (\a4.west);
					\fill[t3] (\a1.center) -- (\a2.center) -- (\a3.center) -- (\a4.center) -- (\a1.center);
				}
			\end{scope}
		\end{tikzpicture}
	}
	\hfill
	\subcaptionbox{subspace $B/C_8$}{
		\begin{tikzpicture}
			\foreach[count=\i, evaluate={\i as \j using int(mod(\i+1,2)+1)}, evaluate={\i as \k using int(mod(\i+2,2)+1)}] \j in {1,2,3,4,1}{
				\node (v\i) at (\i, 1) {$v'_\j$};
				\node (w\i) at (\i,-1) {$v'_\k$};
			}
			\foreach \i in {2,3,4}{
				\draw[halo, t13] (v\i) -- (w\i);
			}
			\draw[halo, t13, ->] (v1) to["$c$"' color=13] coordinate (c1) (w1);
			\draw[halo, t13, ->] (v5) to["$c$"  color=13] coordinate (c2) (w5);
			\draw[t15, halo] ([xshift=2pt]c1) to ([xshift=-2pt]c2);
			\foreach [remember=\i as \j (initially 1)] \i in {2,...,5}{
				\draw[<-, halo, t16] (v\i) to coordinate (g\j) (v\j);
				\draw[<-, halo, t16] (w\i) to coordinate (h\j) (w\j);
				\draw ([yshift=-2pt]g\j) edge[t14, halo] ([yshift=2pt]h\j);
				\coordinate (u\j) at ($(0.5,0)+(\j,0)$);
					\path (u\j) to coordinate (s\j-1) (v\j)
					(u\j) to coordinate (s\j-2) (v\i)
					(u\j) to coordinate (s\j-3) (w\j)
					(u\j) to coordinate (s\j-4) (w\i);
				\scoped[on background layer] {
					\fill[t5] (v\j) -- (s\j-1) -- (g\j) -- (v\j)
						(v\i) -- (s\j-2) -- (g\j) -- (v\i)
						(w\j) -- (s\j-3) -- (h\j) -- (w\j)
						(w\i) -- (s\j-4) -- (h\j) -- (w\i);
					\fill[t6] (u\j) -- (s\j-1) -- (g\j) -- (s\j-2)
						(u\j) -- (s\j-3) -- (h\j) -- (s\j-4);
					\fill[t7] (v\j.center) -- (w\i.center) -- (v\i.center) -- (w\j.center) -- (v\j.center);
				}
			}
		\end{tikzpicture}
	}
	\subcaptionbox{subspace $C/C_8$}{
		\begin{tikzpicture}
			\foreach[count=\j, evaluate={\i as \i using int(mod(\i+1,2)+1)}, evaluate={\l as \l using int(mod(\l+2,2)+1)}] \i/\k/\l in {1/1/3, 2/2/4, 3/3/1, 4/4/2, 1/5/3, 2/6/4, 3/7/1, 4/8/2, 1/1/3} {
				\node (v\j) at (\j, 1) {$v'_\i$};
				\node (w\j) at (\j,-1) {$v'_\l$};
			}
			\foreach \i in {1,...,8}{
				\draw[halo, t17] (v\i) to coordinate[pos=.25] (r\i-1) coordinate[pos=.75] (r\i-2) (w\i);
			}
			\path[t17, ->] (v1) to["$d$"' color=13] coordinate (d1) (w1);
			\path[t17, ->] (v9) to["$d$"  color=13] coordinate (d2) coordinate[pos=.25] (r9-1) coordinate[pos=.75] (r9-2) (w9);
			\draw[t18, halo] ([xshift=2pt]d1) to ([xshift=-2pt]d2);
			\foreach [evaluate=\i as \j using int(\i+1)] \i in {1,...,8}{
				\draw[->, t16, halo] (v\i) to coordinate (g\j) (v\j);
				\draw[->, t16, halo] (w\i) to coordinate (h\j) (w\j);
				\draw ([yshift=-2pt]g\j) edge[t14, halo] coordinate[pos=.125](g\j') coordinate[pos=.875] (h\j') coordinate (s\i) ([yshift=2pt]h\j);
				\scoped[on background layer]{
					\fill[t12] (r\i-1) -- (r\j-2) -- (r\j-1) -- (r\i-2);
					\fill[t8] (r\i-1) -- (r\j-2) -- (h\j') -- (r\i-2) -- (r\j-1) -- (g\j') -- cycle;
					\fill[t10] (r\i-1) -- (g\j') -- (r\j-1) -- (v\j.center) -- (g\j') -- (v\i.center) -- cycle
					(r\i-2) -- (h\j') -- (r\j-2) -- (w\j.center) -- (h\j') -- (w\i.center) -- cycle;
				}
			}
		\end{tikzpicture}
	}
	\caption{The cell structure for $\Cfg(K)$ from \cref{fig:cell structure},
		with vertices replaced by the equivalence classes w.r.t.\ $C_8$ they belong to
		(See \cref{fig:cell structure} for the meaning of the colors).}
	\label{fig:mod-C8:projection}
\end{figure}
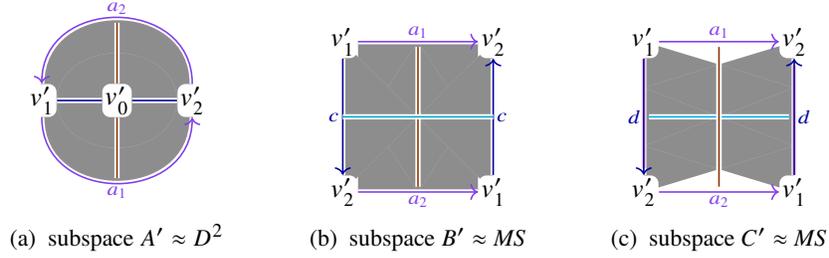
\begin{figure}[tbp]
	\tikzset{x=1cm,y=1cm,baseline=(current bounding box.center)}
	\centering
	\begin{subfigure}{.33\linewidth}\centering
		\begin{tikzpicture}[baseline=(c.base),every edge/.append style={halo}]
			\draw[looseness=1.5]
				(-1,0) 	node (v) {$v'_1$}
				(1,0) 	node (w) {$v'_2$}
				(0,0) 	node (c) {$v'_0$}
				(v.south) 	to[out=-90, in=-90] (w.south) -- (w.north)
				to[out= 90, in= 90] (v.north) -- cycle;
			\draw[looseness=1.5]
				(v) 	edge[out=-90, in=-90, "$a_1$"', ->, t16] coordinate (g) (w)
				(w) 	edge[out= 90, in= 90, "$a_2$"', ->, t16] coordinate (h) (v);
			\draw[t13] (c) edge (v) edge (w);
			\draw[t14] (c) edge coordinate (a) ([yshift=2pt]g)
			               edge coordinate (b) ([yshift=-2pt]h);
			\scoped[on background layer]{
				\fill[t3] (v.east) to[in=180, out=south] (a) to[out=0, in=south] (w.west) to[out=north, in=0] (b) to[out=180, in=north] (v.east);
				\fill[t4] (v.east) to[in=180, out=south] (a) to[out=0, in=south] (w.west) -- (w.south) to[out=south, in=south, looseness=1.5] (v.south) -- cycle
				          (v.east) to[in=180, out=north] (b) to[out=0, in=north] (w.west) -- (w.north) to[out=north, in=north, looseness=1.5] (v.north) -- cycle;
			}
		\end{tikzpicture}
		\caption{subspace $A' \approx D^2$}
	\end{subfigure}%
	\begin{subfigure}{.33\linewidth}\centering
		\begin{tikzpicture}[baseline=(u1.base)]
			\draw[every edge/.append style={halo}]
				(-1, 1) node (v1) {$v'_1$}
				(-1,-1) node (w1) {$v'_2$}
				( 1, 1) node (v2) {$v'_2$}
                ( 1,-1) node (w2) {$v'_1$}
				(0,0) coordinate (u1)
				(v1) edge ["$a_1$" , t16, ->] coordinate (g) (v2)
				(w1) edge ["$a_2$"', t16, ->] coordinate (h) (w2)
				(v1) edge["$c$"', t13, ->] coordinate (t1) (w1)
				(w2) edge["$c$"', t13, ->] coordinate (t2) (v2)
				([yshift=-2pt]g) edge[t14] ([yshift=2pt]h)
				(t1) edge[t15] (t2);
			\coordinate (u) at (0,0);
			\path (u) to coordinate (s-1) (v1)
			      (u) to coordinate (s-2) (v2)
			      (u) to coordinate (s-3) (w1)
			      (u) to coordinate (s-4) (w2);
			\scoped[on background layer] {
				\fill[t5] (v1) -- (s-1) -- (g) -- (v1)
				          (v2) -- (s-2) -- (g) -- (v2)
				          (w1) -- (s-3) -- (h) -- (w1)
				          (w2) -- (s-4) -- (h) -- (w2);
				\fill[t6] (u) -- (s-1) -- (g) -- (s-2)
				          (u) -- (s-3) -- (h) -- (s-4);
				\fill[t7] (v1.center) -- (w2.center) -- (v2.center) -- (w1.center) -- (v1.center);
			}
		\end{tikzpicture}
		\caption{subspace $B' \approx \mathit{MS}$}
	\end{subfigure}%
	\begin{subfigure}{.33\linewidth}\centering
		\begin{tikzpicture}
			\draw[every edge/.append style={halo}]
				(-1, 1) node (v1) {$v'_1$}    (1, 1) node (w1) {$v'_2$}
				(-1,-1) node (w2) {$v'_2$}    (1,-1) node (v2) {$v'_1$}
				(v1) edge[->, "$a_1$", t16] coordinate (g) (w1)
				(w2) edge[->, "$a_2$"', t16] coordinate (h) (v2)
				(v1) edge[t17, ->, "$d$"'] coordinate (t1) coordinate[pos=.25] (r1-1) coordinate[pos=.75] (r1-2) (w2)
				(v2) edge[t17, ->, "$d$"'] coordinate (t2) coordinate[pos=.75] (r2-1) coordinate[pos=.25] (r2-2) (w1)
				([xshift=2pt]t1) edge[t18] ([xshift=-2pt]t2)
				([yshift=-2pt]g) edge[t14, halo] coordinate[pos=.125](g') coordinate[pos=.875] (h') coordinate (s) ([yshift=2pt]h);
			\scoped[on background layer]{
				\fill[t12] (r1-1) -- (r2-2) -- (r2-1) -- (r1-2);
				\fill[t8] (r1-1) -- (r2-2) -- (h') -- (r1-2) -- (r2-1) -- (g') -- cycle;
				\fill[t10] (r1-1) -- (g') -- (r2-1) -- (w1.center) -- (g') -- (v1.center) -- cycle
				           (r1-2) -- (h') -- (r2-2) -- (v2.center) -- (h') -- (w2.center) -- cycle;
			}
		\end{tikzpicture}
		\caption{subspace $C'\approx\mathit{MS}$}
	\end{subfigure}%
	\caption{Cell complex representing $\Cfg(K)/C_8$.}
	\label{fig:cell structure C8}
\end{figure}

\begin{table}[tbp]
	\centering
	$\begin{array}{c >{{}\approx}c c c}
		\toprule
		\multicolumn{2}{c}{\text{subspace}}                   & \beta_\bullet(-, \F_2) & \beta_.(-, \F_3) \\ \midrule
		A'                     & D^2                                           & (1,0,0)                & (1,0,0)          \\
		B', C'                 & \mathit{MS}                                   & (1,0,0)                & (1,0,0)          \\
		A' \cup B', A' \cup C' & \RP^2                                         & (1,1,1)                & (1,0,0)          \\
		B' \cup C'             & \mathit{KB}                                   & (1,2,1)                & (1,1,0)          \\
		A' \cup B' \cup C'     & D^2\cup_{S_1}\mathit{MS}\cup_{S_1}\mathit{MS} & (1,1,2)                & (1,0,0)          \\ \bottomrule
	\end{array}$
	\caption{Expected Betti numbers of various subspaces of $\Cfg(K)/C_8$, according to the cell complex from \cref{fig:cell structure C8}.}
		\label{tab:Betti numbers Cfg/C8}
\end{table}

\begin{figure}[tbp]
		\makebox[\linewidth][c]{
			\subcaptionbox{$\Cfg(K)/C_8$}{
				\includegraphics[scale=0.5, trim={2.5cm 0 5mm 0}, clip]{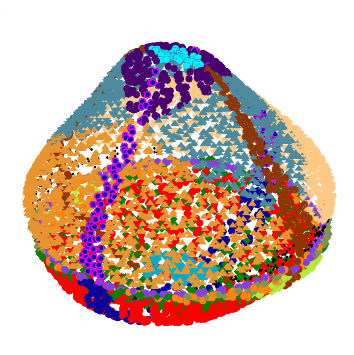}%
				\includegraphics{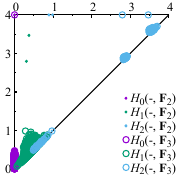}
			}
			\subcaptionbox{$A/C_8$}{
				\includegraphics[scale=0.5, trim={2.5cm 0 05mm 0}, clip]{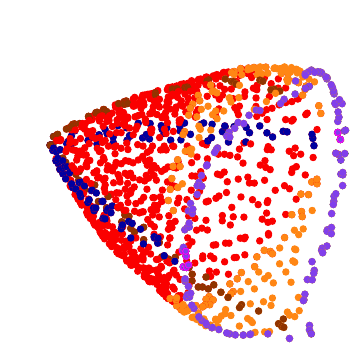}%
				\includegraphics{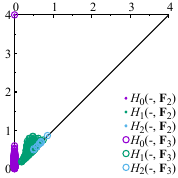}%
			}
		}
		\makebox[\linewidth][c]{
			\subcaptionbox{$B/C_8$\label{fig:B/C8}}{
				\includegraphics[scale=0.5, trim={2.5cm 0 5mm 0}, clip]{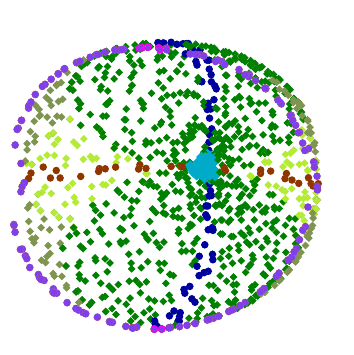}%
				\includegraphics{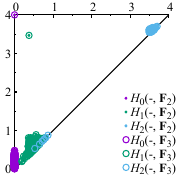}%
			}
			\subcaptionbox{$C/C_8$\label{fig:C/C8}}{
				\includegraphics[scale=0.5, trim={2.5cm 0 5mm 0}, clip]{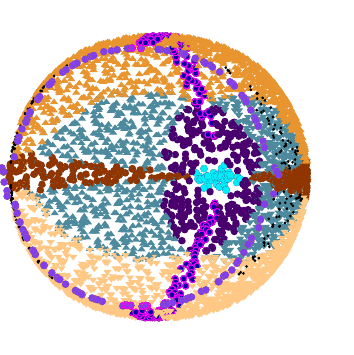}%
				\includegraphics{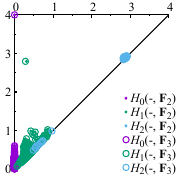}%
			}
		}
		\makebox[\linewidth][c]{
			\subcaptionbox{$(A \cup B)/C_8$}{
				\includegraphics[scale=0.5, trim={2.5cm 0 5mm 0}, clip]{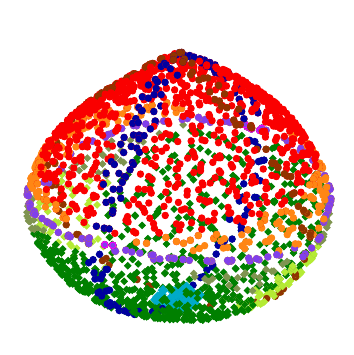}%
				\includegraphics{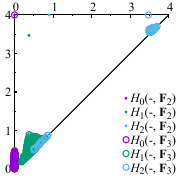}%
			}
			\subcaptionbox{$(A \cup C)/C_8$}{
				\includegraphics[scale=0.5, trim={2.5cm 0 5mm 0}, clip]{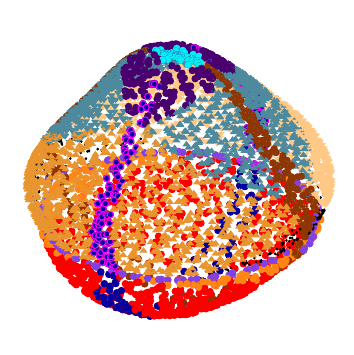}%
				\includegraphics{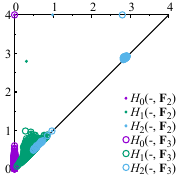}%
			}
		}
		\makebox[\linewidth][c]{
			\subcaptionbox{$(B \cup C)/C_8$}{
				\includegraphics[scale=0.5, trim={2.5cm 0 5mm 0}, clip]{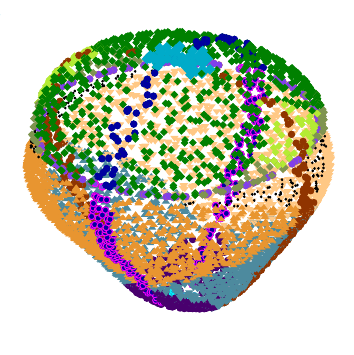}%
				\includegraphics{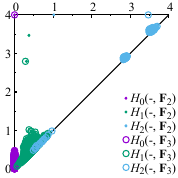}%
			}
			\subcaptionbox*{}[\widthof{%
				\includegraphics[scale=0.5, trim={2.5cm 0 5mm 0}, clip]{new_images/6000.B,C/C8_symmetry/angular2/isomap}%
				\includegraphics{new_images/6000.B,C/C8_symmetry/angular2/intervals}%
			}]{
				\centering
				\includegraphics[scale=0.5]{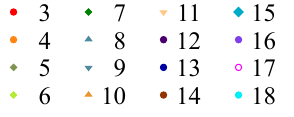}%
			}
		}
	\caption{Isomap projections and persistence diagram (w.r.t. $d_\angle$) of different subspaces of $\Cfg(K)/C_8$.}
	\label{fig:C8}
\end{figure}

\subsection{Plausibility check}
If the model is correct, then both $A' \cup B'$ and $A' \cup C'$ are homeomorphic to a real projective plane $\RP^2 \cong D^2\cup_{S_1}\mathit{MS}$.
Similarly, $B' \cup C'$ is homeomorphic to a Klein bottle $\mathit{KB} \cong \mathit{MS}\cup_{S_1}\mathit{MS}$.
We should thus expect to find the Betti numbers listed in \cref{tab:Betti numbers Cfg/C8}.
In \cref{fig:C8}, we see the persistence diagrams and Isomap projections of these subspaces.
Indeed, we observe the expected Betti numbers.

\begin{figure}[tbp]
	\subcaptionbox{$B/C_8$}{\includegraphics[height=3cm]{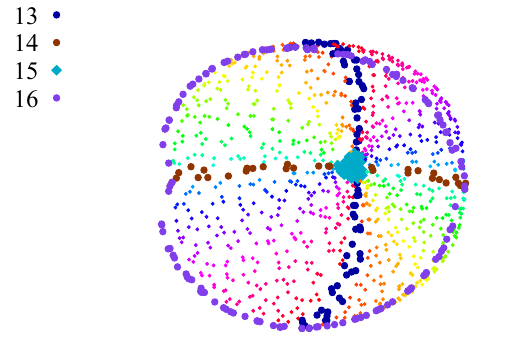}}\hfill
	\subcaptionbox{$C/C_8$}{\includegraphics[height=3cm]{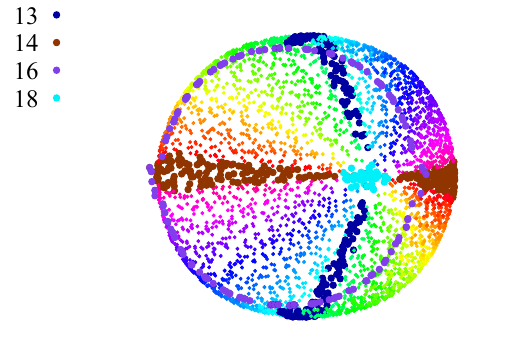}}\hfill
	\subcaptionbox*{}{\includegraphics[height=3cm]{new_images/colorbox}}
	\caption{Color-coded maps $(-) \to S^1$ representing the respective unique most prominent 1-cohomology class in \cref{fig:B/C8,fig:C8}.
		See \cref{sec:circular coordinates} for a description of the method.}
	\label{fig:circular-coordinates/C8}
\end{figure}

To produce further evidence supporting our model,
we have computed maps $B/C_8 \to S^1$ and $C/C_8 \to S^1$ of the presumable Möbius strips $B'$ and $C'$
that correspond to the unique most prominent 1-cohomology class in \cref{fig:B/C8,fig:C/C8}, respective;
see \cref{fig:circular-coordinates/C8}.
We see that the computed maps to $S^1$ wind twice around the common boundary curve ($\Cfg(K, \{16\})/C_8$)
and once about the respective core curves ($\Cfg(K, \{14\})/C_8$ and $\Cfg(K, \{18\})/C_8$)
of $B/C_8$ and $C/C_8$,
which supports our hypothesis.

\section{A cell structure for the quotient space \texorpdfstring{$\Cfg(K)/D_8$}{Cfg(K)/D8}}
\label{sec:Configuration-space-mod-D8}
Finally, we want to understand the unlabeled configuration space $\Cfg(K)/D_8$;
this is the space of all possible configurations
of a chemical molecule with indistinguishable atoms.
As $D_8=C_8\rtimes C_2$, where $C_2$ acts on $C_8$ by inversion.
we have to describe the action of $C_2$ on configurations.

\subsection{Actions}
\label{sec:D8 action-on-configurations}
The generator $s$ of $C_2$ acts on dihedral angles of a $\bar{x} \in \Cfg(K)$ by reversing the order; that is,
\begin{equation}
	\label{eq:C2 action-on-torsion-angles}
	s(\sigma(\bar{x}))_i = \sigma(\bar{x})_{1-i}
\end{equation}
with indices taken modulo $8$.
Analogously, $s$ acts on standard realizations $x^\text{std}$ by
\begin{equation}
	(sx)_i = S x_{1-i}\quad\text{for}\quad S=\begin{psmallmatrix}0 & 1 & 0\\ 1 & 0 & 0\\0 & 0 & -1\end{psmallmatrix};
\end{equation}
i.e., by rotating them by $\pi$ about the axis $(1,1,0)^\mathrm{T}$.
The reindexing and the rotation matrix $S$ are chosen to preserve the Eckart condition \cref{eq:Eckart condition}.
We obtain a Euclidean and an angular quotient metric on $\Cfg(K)/D_8$ by minimizing over actions of $s$ and $t$ analogously to the description in \cref{sec:C8 action-on-configurations}.

\subsection{Cell structure}
\label{sec:cell complex D8}
\begin{figure}[tbp]
	\setlength{\TempLenA}{.25\linewidth}
	\tikzset{x=1cm,y=1cm,baseline=(current bounding box.center)}
	\centering
	\subcaptionbox{subspace $A' \approx D^2$ of $\Cfg(K)/C_8$}[\TempLenA]{
		\begin{tikzpicture}[baseline=(c.base),every edge/.append style={halo}]
			\draw[looseness=1.5]
				(-1,0) 	node (v) {$v'_1$}
				(1,0) 	node (w) {$v'_2$}
				(0,0) 	node (c) {$v'_0$}
				(v.south) 	to[out=-90, in=-90] (w.south) -- (w.north)
				to[out= 90, in= 90] (v.north) -- cycle;
			\draw[looseness=1.5]
				(v) 	edge[out=-90, in=-90, "$a_1$"' near start, ->, t16] coordinate (g) (w)
				(w) 	edge[out= 90, in= 90, "$a_2$"' near start, ->, t16] coordinate (h) (v);
			\draw[t13] (c) edge (v) edge (w);
			\draw[t14] (c) edge coordinate (a) ([yshift=2pt]g)
				edge coordinate (b) ([yshift=-2pt]h);
			\scoped[on background layer]{
				\fill[t3] (v.east) to[in=180, out=south] (a) to[out=0, in=south] (w.west) to[out=north, in=0] (b) to[out=180, in=north] (v.east);
				\fill[t4] (v.east) to[in=180, out=south] (a) to[out=0, in=south] (w.west) -- (w.south) to[out=south, in=south, looseness=1.5] (v.south) -- cycle
				(v.east) to[in=180, out=north] (b) to[out=0, in=north] (w.west) -- (w.north) to[out=north, in=north, looseness=1.5] (v.north) -- cycle;
			}
			\scoped[on background layer]
			\draw [dashed] ([yshift=2em]h) -- ([yshift=-2em]g);
		\end{tikzpicture}
	}\hfill
	\subcaptionbox{subspace $B' \approx \mathit{MS}$ of $\Cfg(K)/C_8$}[\TempLenA]{%
		\begin{tikzpicture}
			\path[every edge/.append style={halo}]
				(-1, 1) node (v1) {$v'_1$}
				( 1, 1) node (v2) {$v'_2$}
				(-1,-1) node (w1) {$v'_2$}
				( 1,-1) node (w2) {$v'_1$}
				( 0, 0) coordinate (u)
				(v1) edge ["$a_1$"  near end,   t16, ->] coordinate (g) (v2)
				(w1) edge ["$a_2$"' near start, t16, ->] coordinate (h) (w2)
				(v1) edge["$c$"' near end, t13, ->] node (t1) {$x''$} (w1)
				(w2) edge["$c$"' near end, t13, ->] node (t2) {$x''$} (v2)
				(g) edge[t14, shorten >=2pt, shorten <=2pt] (h)
				(u) to coordinate (s-1) (v1)
				(u) to coordinate (s-2) (v2)
				(u) to coordinate (s-3) (w1)
				(u) to coordinate (s-4) (w2)
				(t1.east) edge[t15] (t2.west)
				(u) node {$u'$};
			\scoped[on background layer] {
				\path ([yshift=20pt]g) edge[dashed] ([yshift=-20pt]h);
				\fill[t5] (v1) -- (s-1) -- (g) -- (v1)
					(v2) -- (s-2) -- (g) -- (v2)
					(w1) -- (s-3) -- (h) -- (w1)
					(w2) -- (s-4) -- (h) -- (w2);
				\fill[t6] (u) -- (s-1) -- (g) -- (s-2)
					(u) -- (s-3) -- (h) -- (s-4);
				\fill[t7] (v1.center) -- (w2.center) -- (v2.center) -- (w1.center) -- (v1.center);
			}
		\end{tikzpicture}
	}\hfill
	\subcaptionbox{subspace $C'\approx\mathit{MS}$ of $\Cfg(K)/C_8$}[\TempLenA]{
		\begin{tikzpicture}
			\draw[every edge/.append style={halo}]
				(-1, 1) node (v1) {$v'_1$}    (1, 1) node (w1) {$v'_2$}
				(-1,-1) node (w2) {$v'_2$}    (1,-1) node (v2) {$v'_1$}
				(v1) edge[->, "$a_1$"  near end,   t16] coordinate (g) (w1)
				(w2) edge[->, "$a_2$"' near start, t16] coordinate (h) (v2)
				(v1) edge[t17, ->, "$d$"' near end] node (t1) {$t'$} coordinate[pos=.25] (r1-1) coordinate[pos=.75] (r1-2) (w2)
				(v2) edge[t17, ->, "$d$"' near end] node (t2) {$t'$} coordinate[pos=.75] (r2-1) coordinate[pos=.25] (r2-2) (w1)
				(t1.east) edge[t18] (t2.west)
				([yshift=-2pt]g) edge[t14, halo] coordinate[pos=.125](g') coordinate[pos=.875] (h') coordinate (s) ([yshift=2pt]h)
				(s) node {$s'$};
			\scoped[on background layer]{
				\draw ([yshift=20pt]g) edge[dashed] ([yshift=-20pt]h);
				\fill[t12] (r1-1) -- (r2-2) -- (r2-1) -- (r1-2);
				\fill[t8] (r1-1) -- (r2-2) -- (h') -- (r1-2) -- (r2-1) -- (g') -- cycle;
				\fill[t10] (r1-1) -- (g') -- (r2-1) -- (w1.center) -- (g') -- (v1.center) -- cycle
					(r1-2) -- (h') -- (r2-2) -- (v2.center) -- (h') -- (w2.center) -- cycle;
			}
		\end{tikzpicture}
	}
	\\[1em]
	\subcaptionbox{subspace $A''$ of $\Cfg(K)/D_8$}[\TempLenA]{%
		\begin{tikzpicture}
			\path[use as bounding box] (-.2,-1.1) rectangle (2.2,1.1);
			\node (v0) at (0,0) {$v''_0$};
			\node (v) at (2,0) {$v''_1$};
			\coordinate (s1) at (1,1);
			\coordinate (s2) at (1,-1);
			\draw[every edge/.append style={halo}]
				(v0) edge[t14, out=north, in=west] coordinate (r1) (s1) (s1) edge[out=east, in=north, ->, "$a''_2$" , t16] (v)
				(v0) edge[t14, out=south, in=west] coordinate (r2) (s2) (s2) edge[out=east, in=south, ->, "$a''_1$"', t16] (v)
				(v0) edge[t13] (v);
			\scoped[on background layer] {
				\fill[t3]
					(v0.center) to[out=north, in=south west] (r1) to (v.center) to (v0.center)
					(v0.center) to[out=south, in=north west] (r2) to (v.center) to (v0.center);
				\fill[t4]
					(r1) to[out= 30, in=north, in looseness=1.5] (v.center) to (r1)
					(r2) to[out=-30, in=south, in looseness=1.5] (v.center) to (r2)
				;
			}
		\end{tikzpicture}
	}\hfill
	\subcaptionbox{subspace $B''$ of $\Cfg(K)/D_8$}[\TempLenA]{%
		\begin{tikzpicture}
			\path[use as bounding box] (-.2,-1.1) rectangle (2.2,1.1);
			\node (u) at (0,0) {$u''$};
			\node (x) at (1,0) {$x''$};
			\node (v) at (2,0) {$v''_1$};
			\coordinate (s1) at (1,1);
			\coordinate (s2) at (1,-1);
			\draw[every edge/.append style={halo}]
				(u) edge[t14, out=north, in=west] (s1)
				(u) edge[t14, out=south, in=west] (s2)
				(s1) edge[out=east, in=north, ->, "$a''_1$" , t16] coordinate (g) (v)
				(s2) edge[out=east, in=south, ->, "$a''_2$"', t16] coordinate (h) (v)
				(u) edge[t15] (x)
				(x) edge[t13, "$c$"] (v);
			\scoped[on background layer] {
				\path
					(u.north) to[out= 30, in= 150] coordinate(a) (v.north)
					(u.south) to[out=-30, in=-150] coordinate(b) (v.south);
				\fill[t7]
					(u.north) to[out= 30, in=west] (a) to[out=east, in= 150] (v.north) to (v.center)
					(u.south) to[out=-30, in=west] (b) to[out=east, in=-150] (v.south) to (v.center);
				\fill[t6]
					(u.north) to[out=north, in=north west] (g) to (a) to[out=west, in= 30] (u.north)
					(u.south) to[out=south, in=south west] (h) to (b) to[out=west, in=-30] (u.south);
				\fill[t5]
					(a) to (g) to[in=90, out=-30] (v.north) to[out= 150, in=0] (a)
					(b) to (h) to[in=90, out= 30] (v.south) to[out=-150, in=0] (b);
			}
		\end{tikzpicture}
	}\hfill
	\subcaptionbox{subspace $C''$ of $\Cfg(K)/D_8$}[\TempLenA]{%
		\begin{tikzpicture}[x radius=1, y radius=.9]
			\path[use as bounding box] (-1.2,-1.1) rectangle (1.2,1.1);
			\path
				(-1  , 0) node (s) {$s''$}
				(-0.5, 0) node (t) {$t''$}
				( 1  , 0) node (v) {$v''$};
			\draw[t18, halo]
				(s) to (t);
			\draw[t13, halo]
				(t) to coordinate[near start](u) (v);
			\draw[t14, halo]
				(s.north) arc[start angle= 180, end angle= 45] coordinate (a)
				(s.south) arc[start angle=-180, end angle=-45] coordinate (b);
			\draw[t16, halo]
				(a) arc[start angle= 45, end angle=0] coordinate[midway, label={[text=16]above right:\scriptsize$a_1''$}] {}
				(b) arc[start angle=-45, end angle=0] coordinate[midway, label={[text=16]below right:\scriptsize$a_2''$}] {};
			\path
				(s.north) arc[start angle= 180, end angle= 70] coordinate (g)
				(s.south) arc[start angle=-180, end angle=-70] coordinate (h);
			\scoped[on background layer]{
			\fill[t12]
				(s.center) to[out= 90, in= 90] (u)
				(s.center) to[out=-90, in=-90] (u);
			\fill[t8]
				(s.center) -- (s.north) arc[start angle= 180, end angle= 70]  to (u) to[out= 90, in= 90] (s.center)
				(s.center) -- (s.south) arc[start angle=-180, end angle=-70]  to (u) to[out=-90, in=-90] (s.center);
			\fill[t10]
				(u) to (g) to (v.center) to(h) to (u);
			}
		\end{tikzpicture}
	}
	\caption{%
		Upper row: cell $\Cfg(K)/C_8$ from \cref{fig:cell structure C8}. The $D_8$-action takes mirror images across the dashed vertical line.
		With respect to \cref{fig:cell structure C8}, additional nodes and sub-edges have been added.
		Lower row: quotient cell complex structure $\Cfg(K)/D_8$.
		The dashed lines are free edges that are the images of the dashed lines in the upper row.
		The three subspaces are glued along the two lines $a_1''$, $a_2''$,
		which are the image of the edges $a_1'$ and $a_2'$ under the projection.
	}
	\label{fig:cell structure D8}
\end{figure}
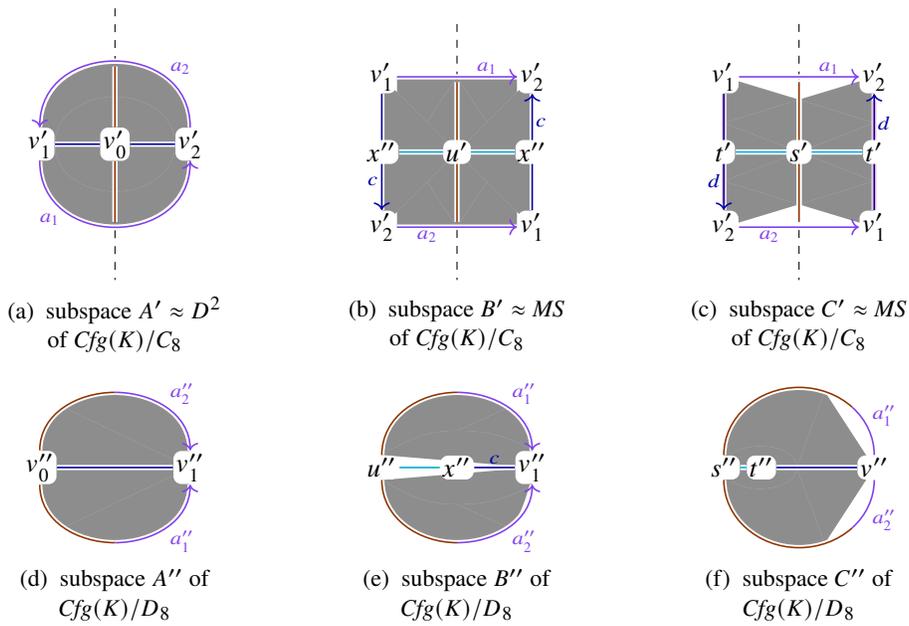

\begin{figure}[tbp]
	\makebox[\linewidth][c]{
		\subcaptionbox{$\Cfg(K)/D_8$}{
			\includegraphics[scale=0.5, trim={2.5cm 0 5mm 0}, clip]{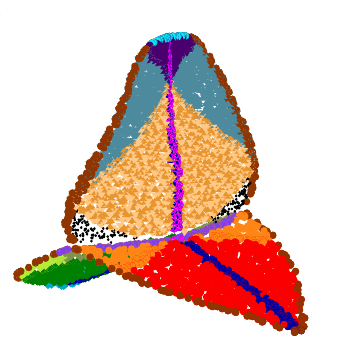}%
			\includegraphics{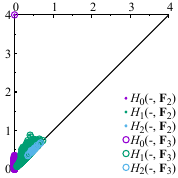}%
		}%
		\subcaptionbox{$A/D_8$}{
			\includegraphics[scale=0.5, trim={2.5cm 0 5mm 0}, clip]{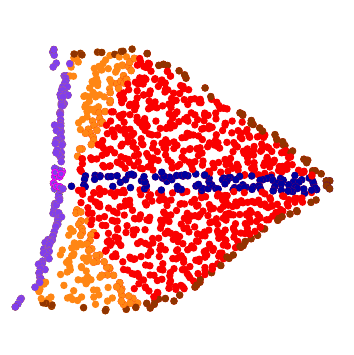}%
			\includegraphics{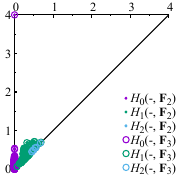}%
		}
	}
	\makebox[\linewidth][c]{
		\subcaptionbox{$B/D_8$}{
			\includegraphics[scale=0.5, trim={2.5cm 0 5mm 0}, clip]{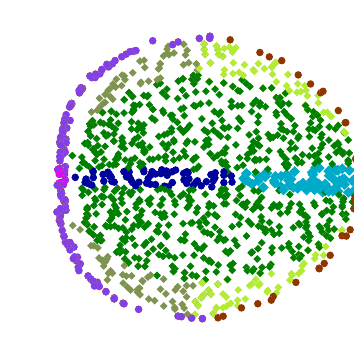}%
			\includegraphics{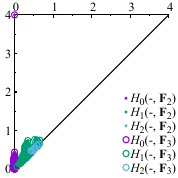}%
		}%
		\subcaptionbox{$C/D_8$}{
			\includegraphics[scale=0.5, trim={2.5cm 0 5mm 0}, clip]{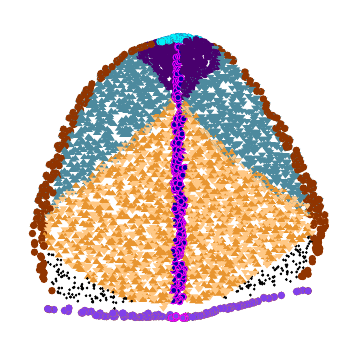}%
			\includegraphics{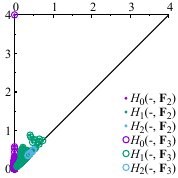}%
		}
	}
	\makebox[\linewidth][c]{
		\subcaptionbox{$(A \cup B)/D_8$}{
			\includegraphics[scale=0.5, trim={2.5cm 0 5mm 0}, clip]{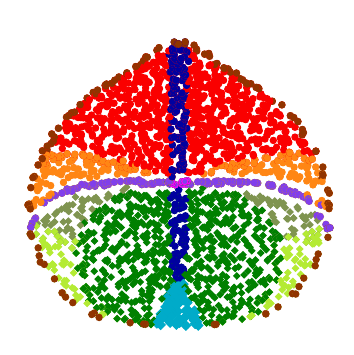}%
			\includegraphics{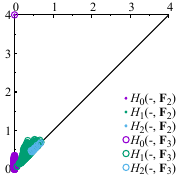}%
		}%
		\subcaptionbox{$(A \cup C)/D_8$}{
			\includegraphics[scale=0.5, trim={2.5cm 0 5mm 0}, clip]{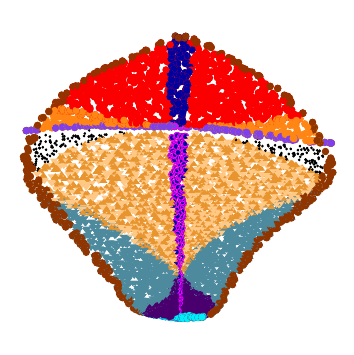}%
			\includegraphics{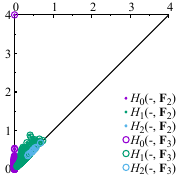}%
		}%
	}
	\makebox[\linewidth][c]{
		\subcaptionbox{$(B \cup C)/D_8$}{
			\includegraphics[scale=0.5, trim={2.5cm 0 5mm 0}, clip]{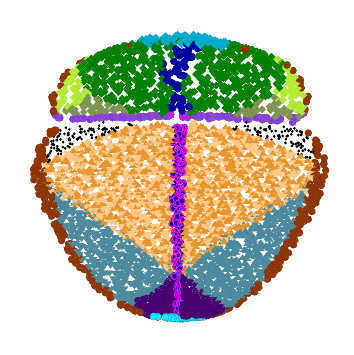}%
			\includegraphics{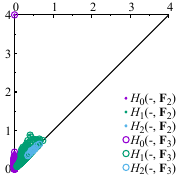}%
		}%
		\subcaptionbox*{}[\widthof{%
			\includegraphics[scale=0.5, trim={2.5cm 0 5mm 0}, clip]{new_images/6000.B,C/C8_symmetry/angular2/isomap}%
			\includegraphics{new_images/6000.B,C/C8_symmetry/angular2/intervals}%
		}]{
			\centering
			\includegraphics[scale=0.5]{new_images/key-horizontal.pdf}%
		}
	}
	\caption{Isomap projections and persistence diagram (w.r.t. $d_\angle$) of different subspaces of $\Cfg(K)/D_8$.}
	\label{fig:D8}
\end{figure}

Recall that the two equivalence classes $v'_1$ and $v'_2$ in $\Cfg(K)/C_8$
correspond to the two torsion angle patterns $(\beta, 0, 0, -\beta, \beta, 0, 0, -\beta)$ for $\beta \gtrless 0$.
Since the $C_2$ action reverses the order of the angles, $v'_1$ and $v'_2$ belong to the same class in $\Cfg(K)/D_8$.
Furthermore, the samples that are both of types 1 and 14 have the dihedral angle pattern $(\alpha, \beta, \alpha, \gamma,\dotsc)$, which is invariant under the action of $C_2$.
Therefore, the quotient map $\Cfg(K)/C_8\twoheadrightarrow \Cfg(K)/D_8$
identifies points as indicated in \cref{fig:cell structure D8}.
The upper row of the figure shows the cell complex representing $\Cfg(K)/C_8$ from \cref{fig:cell structure C8}.
The quotient map $\Cfg(K)/C_8\twoheadrightarrow \Cfg(K)/D_8$ identifies mirror images across the dashed line (type 14).
The lower row contains the image of that cell structure and thus represents a cell complex for $\Cfg(K)/D_8$.
The space $\Cfg(K)/D_8$ thus consists of three sheets glued together along a common line (type 16).
This line is the image of the intersection $A \cap B \cap C$, which is the set of configurations of symmetry type 16.
In particular, $\Cfg(K)/D_8$ is a contractible space.
This is confirmed by the persistence diagrams in \cref{fig:D8}.

\section{Summary and Outlook}
We have investigated the topology of the configuration space
of the cyclic eight-vertex linkage
that corresponds to the chemical cyclooctane molecule,
based on the two metrics $d_\Vert$ and $d_\angle$.
Only the latter actually descends to a metric on the linkage configurations
as is invariant under global translations and rotations
of the configuration.

Using Isomap and persistent homology
we have identified the topology of various subsets of the labeled configuration space $\Cfg(K)$
determined by the symmetries of their elements.
These subsets suggest a cell structure on the space
that is compatible with the quotient maps and hence also descends to
a cell structure of the unlabeled configuration space $\Cfg(K)/D_8$.
The latter space turned out to be a contractible surface consisting of three sheets glued together.

We have not yet considered the configurational energy functional on the conformation space of cyclooctane.
It is an interesting question to study the critical points of this function in light of the cell structure proposed here.

\sloppy
\printbibliography
\end{document}